%% file: main.tex
\documentclass[preprint,12pt]{elsarticle}
\journal{Applied and Computational Harmonic Analysis}

\usepackage{amsmath,amssymb,mathtools,bm}
\usepackage{amsthm}
\usepackage{enumitem}

\usepackage{graphicx}

\usepackage{pgfplots}
\pgfplotsset{compat=1.18}

\usepackage{microtype}

\usepackage{hyperref}
\usepackage[capitalize,nameinlink]{cleveref}
\usepackage{refcount}

\newcommand{\T}{\mathbb{T}}
\newcommand{\R}{\mathbb{R}}
\newcommand{\C}{\mathbb{C}}
\newcommand{\N}{\mathbb{N}}

\newcommand{\e}{\mathrm{e}}
\newcommand{\ii}{\mathrm{i}}

\newtheorem{theorem}{Theorem}[section]
\newtheorem{proposition}[theorem]{Proposition}
\newtheorem{lemma}[theorem]{Lemma}
\newtheorem{corollary}[theorem]{Corollary}

\theoremstyle{definition}

\theoremstyle{remark}
\newtheorem{remark}[theorem]{Remark}

\usepgfplotslibrary{groupplots}
\begin{document}
\begin{frontmatter}

\title{Spectral Conditioning of Weighted Multivariate Vandermonde Operators via Radial Beurling--Selberg Approximation}

\author[addr1]{Joseph Gabet\corref{cor1}}
\ead{joseph.gabet@centralesupelec.fr}
\author[addr1]{Maxime Ferreira Da Costa}
\ead{maxime.ferreira@centralesupelec.fr}
\author[addr2]{Kiryung Lee}
\ead{kiryung@ece.osu.edu}
\cortext[cor1]{Corresponding author.}
\address[addr1]{Laboratory of Signals and Systems, CentraleSup\'elec, Universit\'e Paris--Saclay, CNRS\\3 rue Joliot-Curie, 91190 Gif-sur-Yvette, France}
\address[addr2]{Department of Electrical and Computer Engineering, The Ohio State University\\2015 Neil Avenue, Columbus, OH 43210, USA}

\begin{abstract}
We establish explicit eigenvalue bounds for Gram matrices of weighted multivariate Vandermonde operators with separated nodes on a torus. The squared magnitudes of the weights are samples of a nonnegative radial spectral energy density whose profile has bounded variation, finite energy, and finite weighted variation cost; compact support is not required. The estimate holds in every dimension and depends on the periodic Euclidean separation and the weighted variation integral. For compactly supported profiles, our analysis combines a signed layer decomposition with the ball approximants of Holt and Vaaler; the noncompactly supported case follows by truncation. Combining the resulting radial envelopes with a bandlimited sampling identity based on Poisson summation yields the spectral bounds by eliminating off-diagonal interactions and nonzero diagonal aliases. The result covers nonmonotone profiles and discontinuous radial low-pass filters, and it is complemented by a global elementary upper bound on the Holt--Vaaler cost that exhibits two radial regimes.
\end{abstract}

\begin{keyword}
Multivariate Vandermonde operators \sep spectral conditioning \sep radial weights \sep Beurling--Selberg approximation \sep bounded variation \sep super-resolution
\MSC[2020] 15A12 \sep 15A18 \sep 42B10 \sep 41A30 \sep 94A12
\end{keyword}
\end{frontmatter}

\section{Introduction}

Quantitative spectral bounds for Vandermonde-type operators
play a central role in inverse problems involving point
sources (also called spikes).
We consider a weighted multivariate Vandermonde operator $\bm\Phi_{\bm\tau}: \C^r \to \ell^2(\mathbb{Z}^d; \C)$ given by
\begin{equation}\label{eq:Phi}
[\bm\Phi_{\bm\tau}\bm x]_{\bm k}
= w_{\bm k}
 \sum_{j=1}^r x_j\e^{2\pi\ii\langle\bm k/T,\bm\tau_j\rangle},
\qquad \bm k\in\mathbb Z^d,
\end{equation}
where $\bm\tau=(\bm\tau_1,\ldots,\bm\tau_r)$ denotes the sequence of $r$ nodes in $\R^d$ and $w_{\bm k}$ are complex-valued weights. The sampled frequencies $\bm k/T$, $\bm k\in\mathbb Z^d$, are equispaced. They form the lattice $T^{-1}\mathbb Z^d$, a rescaling of the integer lattice by $1/T$.
The associated Gram matrix depends on the weights only through their squared magnitudes
$|w_{\bm k}|^2$.

Such operators arise in multivariate exponential
fitting~\cite{potts_tasche_multivariate},
direction-of-arrival estimation with sensor
arrays~\cite{stoica_nehorai_music}, and imaging problems
such as super-resolution
microscopy~\cite{ehler_multivariate_microscopy}
and spike
deconvolution~\cite{candes_fernandez_granda_superres,DuvalPeyre2015SparseSpikes,gabet_esprit2025}.
Closely related Fourier-sampling structures also arise
in compressive learning of finite mixture models from
characteristic-function
measurements~\cite{keriven_mixture_sketch,gabet_ferreira_sketch2024}.
In these applications, the smallest singular value of $\bm\Phi_{\bm\tau}$ governs the stability of signal recovery, while its condition number determines the relative sensitivity of the inversion.

Quantitative estimates of these quantities in terms of the node separation and the spectral weights are therefore essential for assessing reconstruction stability, especially for closely spaced nodes and in higher dimensions.

The dependence of spectral conditioning on node separation is closely tied to large-sieve inequalities and extremal approximation in Fourier analysis. In particular, the Beurling--Selberg method constructs bandlimited majorants and minorants of functions; see Vaaler~\cite{vaaler1985}. In dimension one, Moitra~\cite{moitra_vandermonde} exploited this connection to derive conditioning estimates for Vandermonde matrices. Further bounds based on large-sieve techniques were obtained by Aubel and B\"olcskei~\cite{AubelBolcskei2019}, while Li and Liao~\cite{li_liao_restricted_fourier} and Batenkov et al.~\cite{BDGY2020} established lower bounds on the smallest singular value of Fourier matrices with clustered nodes in the super-resolution regime.

Weighted Fourier matrices are particularly relevant to the present work. Kunis and Potts~\cite{KP2007} established stability estimates for weighted multivariate Fourier matrices using decay properties of the associated kernels. Ferreira Da
Costa~\cite{ferreira_weighted_fourier} studied the
conditioning of weighted nonharmonic Fourier matrices
using bandlimited approximation techniques. In related
work, Ferreira Da Costa and
Mitra~\cite{ferreira_mitra_stability} connected
Beurling--Selberg-type approximations to spectral
bounds for sensitivity matrices and to the stability
of the Fisher information matrix.
Hockmann~\cite{hockmann2023AnalysisSparse}
subsequently proved the one-dimensional separation
threshold conjectured in~\cite{ferreira_mitra_stability}
and, using Fisher information bounds, extended the
analysis to ball-shaped frequency sets in higher
dimensions.
Second-order
bandlimited approximations for functions of bounded
variation were developed
in~\cite{ferreira_second_order_beurling}, and
a complementary approach based on discrete variation
was proposed in~\cite{gabet_ferreira2026} for
one-dimensional weighted and first-order confluent
Vandermonde systems. These works provide a natural
starting point for the study of the weighted Vandermonde
operators considered here.

In higher dimensions, several works have studied
unweighted Fourier and Vandermonde operators.
Kunis, Nagel, and
Strotmann~\cite{kunis_nagel_strotmann} established
explicit stability bounds for multivariate Vandermonde
matrices with separated nodes.
For ball-shaped frequency sets, Hockmann and
Kunis~\cite{hockmann2023ShortCommunicationa}
constructed a compactly supported auxiliary
function whose Fourier transform is a
minorant of the ball indicator.
Combined with Poisson summation, this construction
yields local and global Lipschitz stability
estimates for the joint recovery of nodes
and amplitudes.
Li~\cite{li_multivariate2025} studied the smallest
singular values of continuous and discrete nonharmonic
Fourier operators with ball- and cube-shaped frequency
domains. More recently, Diab and
Batenkov~\cite{diab_batenkov2026} investigated the
singular-value scaling of multivariate Vandermonde
matrices with clustered nodes.
Related developments in multidimensional bandlimited
approximation include the work of Carneiro and
Littmann~\cite{carneiro_littmann2014} on extremal
approximations for subordinated radial functions
and that of Gon\c{c}alves, Kelly, and
Madrid~\cite{goncalves_kelly_madrid2015}, who constructed
Gaussian majorants and minorants with box-supported
Fourier transforms, with applications to Hilbert-type
inequalities.

In summary, bandlimited approximation yields weighted conditioning bounds in dimension one, whereas the multivariate results above either concern unweighted ball- or cube-shaped frequency sets or rely on decay properties of the associated kernels~\cite{KP2007}. Radial weights arise naturally when the measurement system is isotropic. The ball-shaped frequency sets of~\cite{hockmann2023ShortCommunicationa,li_multivariate2025} correspond to the simplest radial profile, the indicator of a ball. Nonconstant radial profiles arise, for instance, in scalar incoherent imaging with an ideal circular pupil, whose optical transfer function is radial and nonconstant~\cite{PM2018}; for finite Fourier measurements of point emitters, its squared magnitude defines a compactly supported radial profile of bounded variation. This motivates the study of weighted multivariate Vandermonde operators with radial spectral profiles of bounded variation. For compactly supported profiles, we exploit the radial structure using the ball approximants of Holt and Vaaler~\cite{holtvaaler1996}. Their Fourier transforms are supported in Euclidean balls, and the integrated gap between the upper and lower approximants is explicit. Vaaler~\cite{vaaler2023lattice} used these functions to estimate lattice-point counts in balls. Here, they provide the building blocks for controlling quadratic forms associated with separated exponential sums. A truncation argument then extends the resulting spectral estimate to profiles without compact support.

The main contribution is a spectral estimate for weighted Vandermonde operators whose radial profiles are nonnegative,
of bounded variation, and of finite energy, including nonmonotone and noncompactly
supported profiles whenever their weighted variation cost is finite. In the compactly supported case, the profile is represented as a signed superposition of balls. Its positive and negative variations select the appropriate upper and lower ball approximants. When their Fourier support is matched to the node separation, the approximants eliminate the off-diagonal interactions after Poisson summation. The remaining error is bounded by a single integral of the Holt--Vaaler cost against the radial variation measure. The general case follows by truncation without changing this bound. In dimension one, a sharper variation-based conditioning estimate connects the radial result with the preceding theory of weighted nonharmonic Fourier matrices.

The paper is organized as follows.
Section~\ref{sec:main-results} states the main spectral
bound, develops its low-pass and annular specializations,
and derives an elementary estimate for the Holt--Vaaler
cost. Section~\ref{sec:radial} develops
the radial approximation principle.
Section~\ref{sec:spectral} proves the spectral estimate
and discusses the one-dimensional improvement.
Section~\ref{sec:discussion} discusses perspectives.
Technical details concerning pointwise representatives,
sampling, and the Holt--Vaaler cost, as well as the
proof of the elementary two-regime bound, are collected
in the appendices, together with comparisons
with~\cite{li_multivariate2025}
and~\cite{goncalves_kelly_madrid2015}.

\medskip

\paragraph{Algebraic notation}
Vectors, matrices, and linear operators are written in boldface. For $\bm s,\bm t\in\R^d$, $\langle\bm s,\bm t\rangle$ denotes the standard inner product. For a complex matrix $\bm A$, its conjugate transpose is $\bm A^\ast$, and $\|\bm A\|$ is its operator norm induced by the Euclidean norms. For Hermitian matrices $\bm A,\bm B$ of the same size, $\bm A\preceq\bm B$ means that $\bm B-\bm A$ is positive semidefinite. For vectors, $\|\cdot\|_2$ denotes the Euclidean norm.

\medskip
\paragraph{Asymptotic notation}

Let $t\in\R$ and $a\in\R\cup\{+\infty\}$, and let
$f$ and $g$ be functions of $t$ defined in a suitable
neighborhood of $a$ (or for sufficiently large $t$
when $a=+\infty$). For a stated limit $t\to a$,
we use the standard notations $f=O(g)$, $f\sim g$,
and $f\asymp g$ for an upper bound, asymptotic
equivalence, and two-sided bounds, respectively.
Subscripts, as in $O_d(g)$ or $f\asymp_d g$, indicate
that the implicit constants may depend on the specified
fixed parameters; otherwise, these constants are absolute.
The notation $t\uparrow a$ means that $t$ approaches
$a$ from below. For $a\in\R$, the symbols $t\to a^+$
and $t\to a^-$ denote limits from the right and left,
respectively, while $f(a^+)$ and $f(a^-)$ denote the
corresponding one-sided limits of $f$. In other contexts,
the superscripts $+$ and $-$ have their usual meanings,
as specified by the notation in use.

\medskip
\paragraph{Standing notation}

We denote the sets of nonnegative and positive integers by $\mathbb{N}=\{0,1,2,\ldots\}$ and $\mathbb{N}^\ast=\{1,2,\ldots\}$, respectively. Let $d,r\in\N^\ast$ and $T\in \,]0,+\infty[$.
Let $\T_T^d:=\R^d/(T\mathbb Z^d)$ denote the $d$-dimensional torus of side length $T$.
For $H\in L^1(\R^d;\C)$ and $\bm t\in\R^d$, we use the Fourier convention
\[
\widehat H(\bm t)=\int_{\R^d}H(\bm\xi)\e^{-2\pi\ii\langle\bm\xi,\bm t\rangle}\,\mathrm d\bm\xi.
\]
For $\bm z\in\R^d$ and $\eta\in \,]0,+\infty[$, write
$B_2(\bm z,\eta):=\{\bm t\in\R^d:\|\bm t-\bm z\|_2<\eta\}$
and $\overline B_2(\bm z,\eta)$ for its closure.
The surface area of the unit sphere is
$\omega_{d-1}=2\pi^{d/2}/\Gamma(d/2)$; in particular, $\omega_0=2$.

\section{Main results}\label{sec:main-results}
\subsection{Radially weighted Vandermonde operators}

We consider a special form of the weighted multivariate Vandermonde operator defined in Equation~\eqref{eq:Phi}, where
the weights $w_{\bm k}$ satisfy the radial constraint
\begin{equation}\label{eq:w}
    |w_{\bm k}|^2=\varphi(\|\bm k/T\|_2)
\end{equation}
for a nonnegative function $\varphi$ of bounded variation on $]0,+\infty[$
such that $\varphi(\|\cdot\|_2)$ is integrable on $\R^d$. We refer to $\varphi$ as the \emph{radial profile} of the Vandermonde operator.
Since the pointwise values of the radial profile matter for discrete sampling,
we fix a representative on $[0,+\infty[$ by setting
\begin{equation}\label{eq:representative}
\varphi(\rho)
=\tfrac12\bigl(\varphi(\rho^-)+\varphi(\rho^+)\bigr)
\quad(\rho>0),
\qquad
\varphi(0)=\varphi(0^+).
\end{equation}
Thus, at a jump located at a radius $\rho>0$, we assign the average of the two one-sided limits, while at the origin we use the right-hand limit. This choice fixes the pointwise values of the radial profile without changing its integral or its distributional derivative $D\varphi$. In the sequel, we denote by $|D\varphi|$ its total-variation measure.

We denote by $E$ the integral over $\R^d$ of the radial function associated with the profile,
\begin{equation}\label{eq:energy}
E:=\int_{\R^d}\varphi(\|\bm\xi\|_2)\,\mathrm d\bm\xi
=\omega_{d-1}\int_0^{+\infty}\varphi(\rho)\rho^{d-1}\,\mathrm d\rho,
\end{equation}
so that $T^dE$ serves as a proxy for the energy $\sum_{\bm k}|w_{\bm k}|^2$ of the weight sequence (see Remark~\ref{rem:energy}). We consider only nontrivial cases where $E>0$.

Write $\nu=(d-2)/2$, and let $J_\nu$ denote the Bessel function of the first kind of order $\nu$. The Holt--Vaaler cost $u_\nu$ is defined by
\begin{equation}\label{eq:u}
u_\nu(\rho,\delta)
=\frac{\delta^{-1}\rho^{d-1}}
 {1-\frac{\pi(d-1)}2\displaystyle\int_{\pi\delta\rho}^{+\infty}
       \frac{J_\nu(t)J_{\nu+1}(t)}t\,\mathrm dt}, \qquad \rho,\delta\in\,]0,+\infty[,
\end{equation}
so that $\omega_{d-1}u_\nu(\rho,\delta)$ is the integrated gap between the Holt--Vaaler majorant and minorant of the indicator of a ball of radius $\rho$ with bandwidth $\delta$; see Lemma~\ref{lem:ball-pair} and~\cite{vaaler2023lattice}.
The cost~\eqref{eq:u} is positive and has a continuous extension to $\rho=0$; see \ref{app:cost}. In dimension one, Equation~\eqref{eq:u} reduces to $u_{-1/2}(\rho,\delta)=\delta^{-1}$.

We introduce the weighted variation cost at scale $\delta$, which determines how the
individual ball approximation costs accumulate (see
Section~\ref{sec:radial}):
\begin{equation}\label{eq:W}
W_\delta:=\omega_{d-1}\int_{]0,+\infty[}
u_\nu(\rho,\delta)\,\mathrm d|D\varphi|(\rho).
\end{equation}
Like $E$, the cost $W_\delta$ is determined by the continuous profile, not by its samples alone.

Given nodes $\bm \tau = (\bm\tau_1,\ldots,\bm\tau_r)\in (\T_T^d)^r$, we define the \emph{minimal separation} as
\begin{equation}\label{eq:sep}
\Delta:=\min\left\{T,\ \min_{\substack{j,\ell\in\{1,\ldots,r\}\\j\ne\ell}}\min_{\bm q\in\mathbb Z^d}
 \|\bm\tau_j-\bm\tau_\ell+T\bm q\|_2\right\},
\end{equation}
where differences are taken between any representatives in $\R^d$. The inner minimum over $\bm q$ does not depend on this choice. For $r=1$, the minimum over the empty set of distinct pairs is
$+\infty$, so that $\Delta=T$. We assume that the nodes are pairwise distinct in the torus, so that $0<\Delta\leq T$.

The following theorem shows that the spectrum of the Gram matrix associated with the operator $\bm\Phi_{\bm\tau}$ defined in~\eqref{eq:Phi} is controlled by the continuous energy $T^dE$, up to an error governed by $W_\Delta$.

\begin{theorem}[Radial conditioning bound]\label{th:main}
For an integrable radial profile $\varphi$ of bounded variation such that $W_\Delta<+\infty$, the radially weighted Vandermonde operator with profile $\varphi$ defined by \eqref{eq:Phi} and \eqref{eq:w} is bounded and satisfies
\begin{equation}\label{eq:gram}
T^d(E-W_\Delta)\bm I_r
\preceq\bm\Phi_{\bm\tau}^\ast\bm\Phi_{\bm\tau}
\preceq T^d(E+W_\Delta)\bm I_r,
\end{equation}
where $E$ and $W_\Delta$ are respectively defined in \eqref{eq:energy} and \eqref{eq:W}.
\end{theorem}

The detailed proof is given in Section~\ref{sec:spectral}.
In particular, Theorem~\ref{th:main} provides uniform control of the deviation of the normalized Gram matrix from the identity:
\begin{equation}\label{eq:normalized}
\left\|\frac{\bm\Phi_{\bm\tau}^\ast\bm\Phi_{\bm\tau}}
{T^dE}-\bm I_r\right\|
\leq \frac{W_\Delta}{E}.
\end{equation}
Consequently, the relative error $\varepsilon:=W_\Delta/E$ directly controls the conditioning of the operator. Whenever $\varepsilon<1$, the operator has full column rank and its spectral condition number satisfies
\begin{equation}\label{eq:kappa}
\kappa_2(\bm\Phi_{\bm\tau})
=\sqrt{\frac{
\lambda_{\max}(\bm\Phi_{\bm\tau}^\ast\bm\Phi_{\bm\tau})
}{
\lambda_{\min}(\bm\Phi_{\bm\tau}^\ast\bm\Phi_{\bm\tau})
}}
\leq\sqrt{\frac{1+\varepsilon}{1-\varepsilon}}.
\end{equation}
Thus, $W_\Delta<E$ provides an explicit sufficient condition for stable inversion. %
The comparison in Appendix~\ref{app:comparison-gkm} with the box-supported Gaussian approximants of Gon\c{c}alves, Kelly, and Madrid, for the noncompact Gaussian profile, illustrates the trade-off between the generality of our approach and the sharper spectral bounds obtained from Gaussian-specific approximants in that example.

\begin{remark}[Continuous and sampled energy]\label{rem:energy}
Every diagonal entry of the Gram matrix equals the sampled
energy
\(
S_T:=\sum_{\bm k\in\mathbb Z^d}
\varphi(\|\bm k/T\|_2),
\)
which need not equal $T^dE$. By \eqref{eq:normalized},
\(
\left|\frac{S_T}{T^dE}-1\right|
\leq\frac{W_\Delta}{E}.
\)
Thus the bound controls both the diagonal sampling error
and the off-diagonal correlations. For $r=1$, it remains
a nontrivial estimate of the sampled energy.
\end{remark}

\begin{remark}[Derivatives with respect to the node locations]
Theorem~\ref{th:main} also yields bounds for the operator formed by the spatial gradients of the Vandermonde columns. Indeed, collecting these gradients gives the vector-valued operator
\(
[\bm{\tilde\Phi}_{\bm\tau}]_{\bm k,j}
=
2\pi\ii\,\frac{\bm k}{T}
w_{\bm k}
\e^{2\pi\ii\langle\bm k/T,\bm\tau_j\rangle}.
\)
Its Gram matrix coincides with that of a radially weighted
Vandermonde operator with the modified profile
\(
\tilde\varphi(\rho)=4\pi^2\rho^2\varphi(\rho).
\)
Consequently, Theorem~\ref{th:main} applies directly to
$\bm{\tilde\Phi}_{\bm\tau}^\ast \bm{\tilde\Phi}_{\bm\tau}$,
provided that $\tilde\varphi$ satisfies the required assumptions.
These bounds quantify the conditioning of the columnwise gradient
operator, which naturally arises in the study of position-dependent
Jacobians for spike deconvolution~\cite{gabet_esprit2025}.
This operator differs from the full Jacobian with respect to all node coordinates, whose Gram matrix also contains cross terms between coordinates.
\end{remark}

In dimension one, the classical Beurling--Selberg approximation yields equal upper and lower integrated errors, leading to a factor-of-two improvement over the direct specialization of Theorem~\ref{th:main}.

The resulting conditioning bound coincides, after normalization, with the one-dimensional result previously established by Ferreira Da Costa~\cite[Theorem~1, Eqs.~(6a)--(6b)]{ferreira_weighted_fourier}. Our formulation additionally accommodates noncompact spectral profiles, beyond the finitely many samples considered there.
We state it here to make explicit the one-dimensional specialization of our general framework.

\begin{proposition}\label{prop:one-dimensional}
If $d=1$, then under the assumptions of
Theorem~\ref{th:main},
\begin{equation}\label{eq:one-dimensional}
\left\|
\frac{\bm\Phi_{\bm\tau}^\ast\bm\Phi_{\bm\tau}}{TE}
-\bm I_r
\right\|
\leq
\frac{|D\varphi|(]0,+\infty[)}{E\Delta}.
\end{equation}
\end{proposition}

The detailed proof is given in Section~\ref{sec:one-dimensional-proof}.

\subsection{Ideal radial low-pass and annular band-pass filters}

We now illustrate Theorem~\ref{th:main} with two classical
radial spectral profiles: a low-pass filter supported on a
ball and a band-pass filter supported on an annulus.
In both cases, the variation measure is explicitly known,
leading to concrete conditioning estimates.

For a real radius $R\in\,]0,+\infty[$, write $\chi_R$ for
the indicator of $B_2(0,R)$ with value $1/2$ on its boundary
sphere.

The ideal radial low-pass filter provides a first application
of Theorem~\ref{th:main}. Since its radial profile has a
single jump at the cutoff frequency $R$, the variation integral
reduces to a single Holt--Vaaler cost. If the bandwidth $B>0$
is described by a diameter, then $R=B/2$. This yields the
following explicit estimate.

\begin{corollary}[Radial low-pass]\label{cor:ball}
Suppose that $\varphi(\|\bm\xi\|_2)=\chi_R(\bm\xi)$.
Then $E=\omega_{d-1}R^d/d$ and
\begin{equation}\label{eq:ball}
\left\|\frac{\bm\Phi_{\bm\tau}^\ast\bm\Phi_{\bm\tau}}{T^dE}-\bm I_r\right\|
\leq \frac{d\,u_\nu(R,\Delta)}{R^d}.
\end{equation}
The same bound holds for the open or closed ball indicator,
including when the sampling lattice meets its boundary.
\end{corollary}

\begin{proof}
The radial variation measure is the unit point mass at $R$,
giving \eqref{eq:ball}. The ball approximants are continuous.
Their pointwise inequalities on the interior and exterior
therefore imply that the upper approximant is at least $1$
on the sphere and the lower approximant is at most $0$ there.
The same envelopes sandwich either boundary convention,
with unchanged integrals.
\end{proof}

The bound in \eqref{eq:ball} depends on the node configuration
only through the periodic Euclidean separation $\Delta$
defined in \eqref{eq:sep}, with no explicit dependence on
the number $r$ of nodes. In particular, whenever the
right-hand side is strictly smaller than one, the Gram
matrix is positive definite and the operator has full
column rank. However, a bound greater than or equal
to one does not imply rank deficiency.

The behavior of the estimate becomes particularly simple
when the product $R\Delta$ is large. Indeed, the asymptotic
formulas in \ref{app:cost} yield, for every fixed
dimension $d\in\N^\ast$,
\[
\frac{d\,u_\nu(R,\Delta)}{R^d}
\sim \frac{d}{R\Delta},
\qquad R\Delta\to+\infty.
\]
Thus, in fixed dimension, the relative spectral error is
bounded by a quantity that decays asymptotically like
$(R\Delta)^{-1}$. This equivalence is not asserted uniformly
in the dimension and does not, by itself, describe the
simultaneous regime $d\to+\infty$ with $R\Delta\asymp d$.

In dimension three, the right-hand side of \eqref{eq:ball}
admits the exact expression
\begin{equation}\label{eq:three}
\frac{3}{R\Delta}
\left[
1-\left(\frac{\sin(\pi R\Delta)}{\pi R\Delta}\right)^2
\right]^{-1}.
\end{equation}
The correction factor in \eqref{eq:three} tends to one as $R\Delta\to+\infty$,
consistent with the preceding asymptotic equivalence.

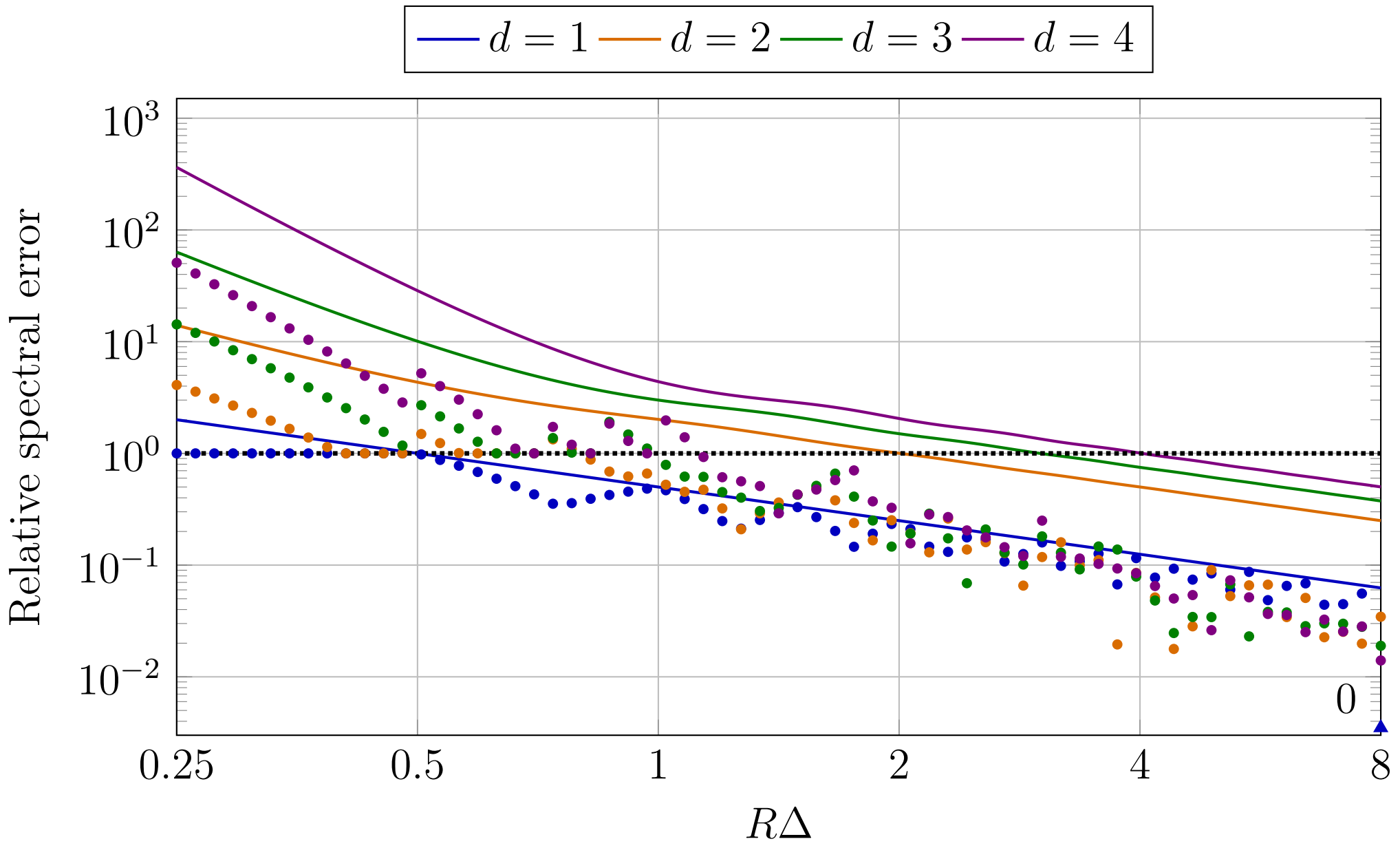
\begin{figure}[htbp]
\centering
\input{figures/radial_lowpass.tex}
\caption{Low-pass spectral bounds versus $R\Delta$
in dimensions $d=1,2,3,4$. Solid curves represent
the sharper bound $1/(2R\Delta)$ from
Proposition~\ref{prop:one-dimensional} when $d=1$,
and the integrated-gap bound
$d\,u_{(d-2)/2}(R,\Delta)/R^d$ from
Corollary~\ref{cor:ball} when $d\geq2$.
Points show the normalized operator-norm deviations
in \eqref{eq:normalized} for the periodic grids
$\{0,\Delta\}^d$, with $R=1$ and $T=2\Delta$.
Theoretical curves and empirical points both cover
$R\Delta\in[0.25,8]$, including the endpoints.
The zero deviation for $d=1$ at $R\Delta=8$
is shown by a triangle labeled $0$ at the bottom
of the logarithmic axis. The dotted horizontal
line marks the sufficient threshold $1$ for a
positive lower spectral bound. The figure-generation code is available in~\cite{GabetRadialVandermondeFigures2026}.}
\label{fig:radial}
\end{figure}

Figure~\ref{fig:radial} illustrates the dependence of \eqref{eq:ball} on the dimensionless separation $R\Delta$ using a simple family of periodic grids. The crossing of the level $1$ marks a sufficient condition for positive definiteness, not a necessary threshold for the actual spectrum. A quantitative comparison with Li's bound on the same ball-sampling Gram matrix is given in \ref{app:comparisons}.

\medskip

The preceding example concerns a nonincreasing radial
profile. However, Theorem~\ref{th:main} does not require
monotonicity. To illustrate this point, consider a radial
band-pass filter whose spectral support is an annulus.

\begin{corollary}[Radial annulus]\label{cor:annulus}
Let $a,b\in\R$ satisfy $0<a<b$, and let the spectral
energy density satisfy $\varphi(\|\bm\xi\|_2)=\chi_b(\bm\xi)-\chi_a(\bm\xi)$. Then
\begin{equation}\label{eq:annulus}
\left\|\frac{\bm\Phi_{\bm\tau}^\ast\bm\Phi_{\bm\tau}}{T^dE}-\bm I_r\right\|
\leq
\frac{d\,[u_\nu(a,\Delta)+u_\nu(b,\Delta)]}{b^d-a^d}.
\end{equation}
\end{corollary}

\begin{proof}
Here $-D\varphi$ is the unit point mass at $b$ minus
the unit point mass at $a$. Its total variation puts
unit mass at both radii, while
$E=\omega_{d-1}(b^d-a^d)/d$.
\end{proof}

Unlike the low-pass profile, the annular profile has two
jumps, of opposite signs, at the radii $a$ and $b$.
Since the error term involves the total variation of
the radial profile, both jumps contribute positively
to the numerator of \eqref{eq:annulus}.

On the other hand, the continuous spectral energy is
proportional to $b^d-a^d$. As the annulus becomes
arbitrarily thin, this energy tends to zero while
the two Holt--Vaaler costs remain positive.
Consequently, the resulting relative error bound can
become large even when the outer radius $b$ is large.
This illustrates the role of total variation in the
estimate, without implying that the bound is sharp
for such spectral profiles.

\begin{remark}[Dimension one]\label{rem:one-dimensional-filters}
For $d=1$, Proposition~\ref{prop:one-dimensional} improves both corollaries. With the boundary values fixed by \eqref{eq:representative}, the low-pass profile has $E=2R$ and $|D\varphi|(]0,+\infty[)=1$, while the annular profile has $E=2(b-a)$ and $|D\varphi|(]0,+\infty[)=2$. Hence
\begin{equation}\label{eq:one-dimensional-filters}
\begin{aligned}
\left\|\frac{\bm\Phi_{\bm\tau}^\ast\bm\Phi_{\bm\tau}}{TE}-\bm I_r\right\|
& \leq\frac{1}{2R\Delta}
\quad\text{(low-pass)},\\
\left\|\frac{\bm\Phi_{\bm\tau}^\ast\bm\Phi_{\bm\tau}}{TE}-\bm I_r\right\|
& \leq\frac{1}{(b-a)\Delta}
\quad\text{(annulus)}.
\end{aligned}
\end{equation}
\end{remark}

\subsection{An elementary bound with two regimes}

The preceding conditioning estimates involve the Holt--Vaaler cost $u_\nu$ defined in \eqref{eq:u}, expressed through an integral
of Bessel functions. We now derive an elementary upper
bound that makes its dependence on the radius, the
dimension, and the parameter $\delta$ explicit.
The resulting estimate reveals two distinct radial
regimes and yields a simple sufficient condition for
a positive lower spectral bound, hence for positive definiteness
of the Gram matrix, in the low-pass case.

\begin{proposition}[Elementary bound for the Holt--Vaaler cost]\label{prop:elementary}
Let $d\in\N^\ast$ satisfy $d\geq2$. Then, for every $\rho\geq0$ and every $\delta>0$,
$u_{(d-2)/2}(\rho,\delta)\leq u_{\rm simp}(\rho,\delta)$, where
\begin{equation}\label{eq:elementary}
u_{\rm simp}(\rho,\delta):=
\begin{cases}
c_d \delta^{-d} +\bigl(2\delta^{-1}\rho_\delta^{d-1}-c_d \delta^{-d} \bigr)\rho/\rho_\delta,
&0\leq\rho\leq\rho_\delta,\\[2mm]
2\delta^{-1}\rho^{d-1},&\rho\geq\rho_\delta,
\end{cases}
\end{equation}
where
\[
\rho_\delta:=\frac{d-1}{\pi\delta},\qquad
c_d:=\frac{d\,2^{d-1}\Gamma(d/2)^2}{\pi^d}.
\]
\end{proposition}

\medskip
The proof of Proposition~\ref{prop:elementary} is given in \ref{app:elementary}.
\medskip

\begin{figure}[htbp]
\centering
\input{figures/elementary_cost.tex}
\caption{The exact Holt--Vaaler cost and the
elementary majorant from
Proposition~\ref{prop:elementary}
for $d=8$ and $\delta=1$.
Top: the exact cost $u_3(\rho,1)$, its elementary
majorant, the large-radius asymptotic $\rho^7$,
and the origin value $u_3(0,1)$.
Bottom: the ratio $u_{\rm simp}(\rho,1)/u_3(\rho,1)$ of the elementary majorant~\eqref{eq:elementary} to the exact cost~\eqref{eq:u}. The vertical dotted line
marks $\rho_\delta=7/\pi$. The figure-generation code is available in~\cite{GabetRadialVandermondeFigures2026}.}
\label{fig:elementary}
\end{figure}
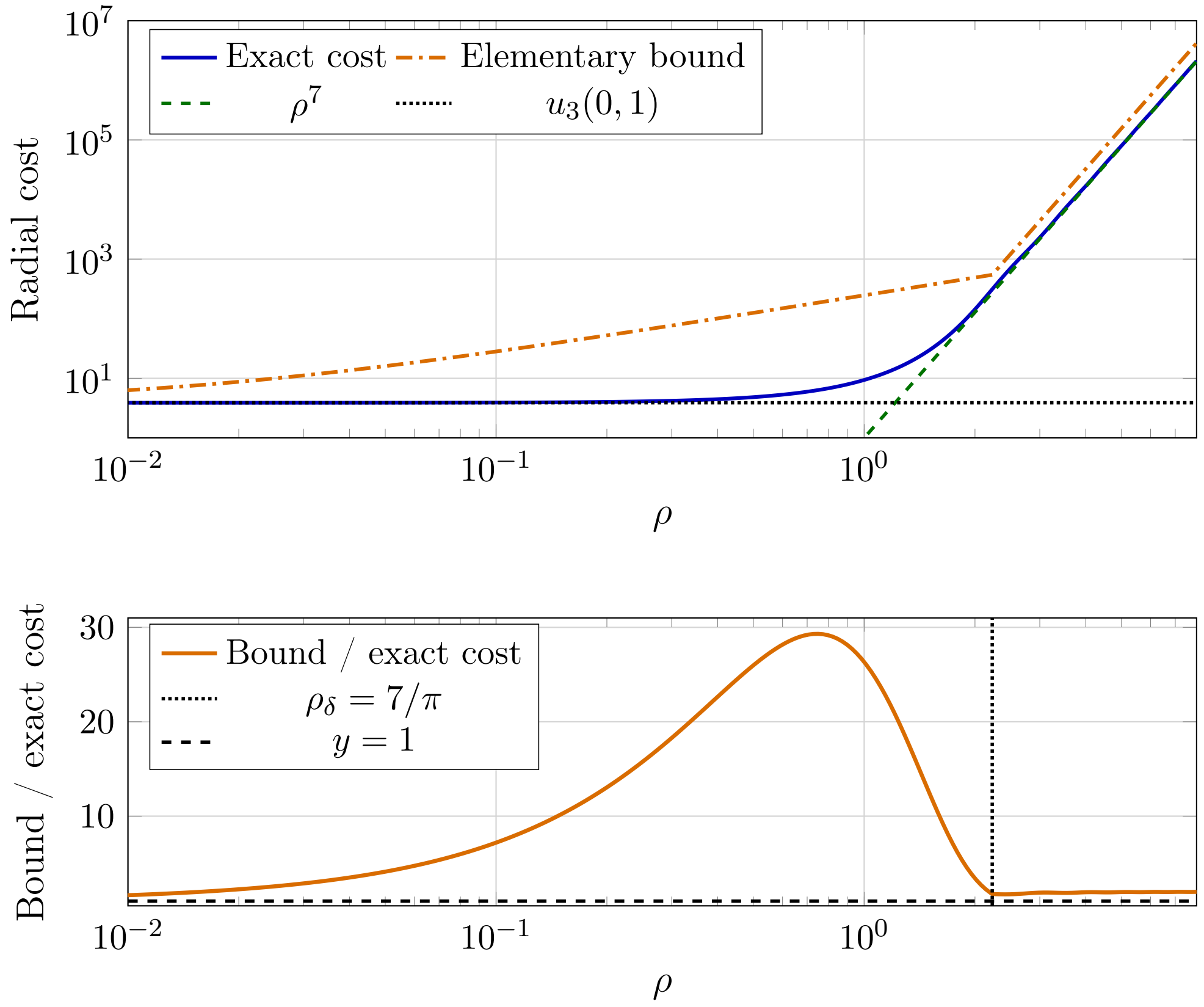

The threshold $\rho_\delta=(d-1)/(\pi\delta)$ separates the two regimes of the elementary majorant in Proposition~\ref{prop:elementary}: an affine bound near the origin and a large-radius bound of order $\delta^{-1}\rho^{d-1}$. Figure~\ref{fig:elementary} illustrates these regimes and the ratio between the majorant and the exact Holt--Vaaler cost. Although the elementary bound captures the correct large-radius order, it may be conservative at intermediate radii. Its validity follows from Proposition~\ref{prop:elementary}, independently of the numerical illustration.

Replacing $u_\nu$ by $u_{\rm simp}$ in \eqref{eq:W} therefore gives a fully elementary sufficient estimate for any of the radial profiles in Theorem~\ref{th:main}. As an example, for the low-pass case, in the \emph{large-separation regime} where $R\Delta \geq \frac{d-1}{\pi}$, Corollary~\ref{cor:ball} yields
\begin{equation}\label{eq:elementary-ball}
\left\|\frac{\bm\Phi_{\bm\tau}^\ast\bm\Phi_{\bm\tau}}{T^dE}-\bm I_r\right\|
\leq \frac{2d}{R\Delta}.
\end{equation}
In particular, the explicit estimate~\eqref{eq:elementary-ball} shows that $R\Delta>2d$ guarantees a positive lower spectral bound. This linear dependence on dimension is consistent with the necessary dimensional growth of universal separation thresholds for interpolation with Euclidean-ball spectra established by Olevskii and Ulanovskii~\cite[Corollary~5.1]{OU2012}. This comparison concerns the dimensional scaling, not a necessary separation threshold for the finite matrices considered here.

For the unweighted ball-sampling model,
Hockmann and
Kunis~\cite[Theorem~3.3]{hockmann2023ShortCommunicationa}
also obtain explicit Lipschitz stability
estimates at a dimensionless separation
scale of $2d$.
Their result concerns joint node-and-amplitude recovery rather than a two-sided bound
on the Gram spectrum.

\section{A radial Beurling--Selberg principle}\label{sec:radial}
\subsection{Ball approximants and the Holt--Vaaler cost}

Our construction of bounds for radial functions relies on bandlimited approximations of indicators of Euclidean balls. The following lemma, due to Holt and Vaaler~\cite[Theorem~3]{holtvaaler1996}, provides a minorant and a majorant for each ball, with prescribed Fourier support and an explicit integrated gap. These approximants serve as building blocks for treating more general radial functions by superposition. We use the normalization recalled in \cite[Lemma~3.1 and equations~(40), (52)--(53)]{vaaler2023lattice}.

\begin{lemma}[Holt--Vaaler ball pair]\label{lem:ball-pair}
Let $d\in\N^\ast$ and $\rho,\delta\in \,]0,+\infty[$. There exist real-valued radial functions $F_\delta(\cdot;\rho),G_\delta(\cdot;\rho)\in C(\R^d;\R)\cap L^1(\R^d;\R)$ such that
\begin{equation}\label{eq:ball-pair}
F_\delta(\bm\xi;\rho)\leq\chi_\rho(\bm\xi)\leq G_\delta(\bm\xi;\rho),
\qquad
\operatorname{supp}\widehat F_\delta,\operatorname{supp}\widehat G_\delta
\subseteq\overline B_2(0,\delta),
\end{equation}
and
\begin{equation}\label{eq:ball-gap}
\int_{\R^d}(G_\delta-F_\delta)(\bm\xi;\rho)\,\mathrm d\bm\xi
=\omega_{d-1}u_\nu(\rho,\delta).
\end{equation}
The explicit ball pair can be chosen measurably in $\rho$.
\end{lemma}

The identity \eqref{eq:ball-gap} gives the exact integrated gap of the Holt--Vaaler pair. This identity will allow us to control the approximation error when superposing ball approximants.

\subsection{Signed radial layer decomposition}
In this section we assume, in addition, that $\varphi$ has
compact support. The resulting approximation pair is used to
prove the compact-support case of Theorem~\ref{th:main}; its
extension to noncompact profiles is obtained separately by
truncation in Section~\ref{sec:spectral}.

For every $\bm\xi\in\R^d$, set $f(\bm\xi):=\varphi(\|\bm\xi\|_2)$, and let $\mu:=-D\varphi$ be a finite signed Borel measure on $]0,+\infty[$. Then
\begin{equation}\label{eq:layers}
f(\bm\xi)=\int_{]0,+\infty[}\chi_\rho(\bm\xi)\,\mathrm d\mu(\rho).
\end{equation}
Indeed, for a real number $s\in \,]0,+\infty[$, the mass strictly above $s$ is $\varphi(s^+)$, and half the mass at $s$ is $[\varphi(s^-)-\varphi(s^+)]/2$. At $s=0$, the total mass is $\varphi(0^+)$. Thus \eqref{eq:representative} gives the identity at every point, not merely almost everywhere.

Write $\mu=\mu^+-\mu^-$ for its Jordan decomposition, where $\mu^+$ and $\mu^-$ are mutually singular positive measures. Negative layers require the opposite envelope. Using the upper ball approximant for both signs would reverse the desired inequality on those layers. For a fixed real number $\delta\in \,]0,+\infty[$, define
\begin{equation}\label{eq:envelopes}
\begin{aligned}
M_\delta^+f&=\int_{]0,+\infty[} G_\delta(\cdot;\rho)\,\mathrm d\mu^+(\rho)
             -\int_{]0,+\infty[} F_\delta(\cdot;\rho)\,\mathrm d\mu^-(\rho),\\
M_\delta^-f&=\int_{]0,+\infty[} F_\delta(\cdot;\rho)\,\mathrm d\mu^+(\rho)
             -\int_{]0,+\infty[} G_\delta(\cdot;\rho)\,\mathrm d\mu^-(\rho).
\end{aligned}
\end{equation}

\begin{proposition}[Radial approximation pair]\label{prop:radial}
The functions in \eqref{eq:envelopes} have radial representatives in $C(\R^d;\R)\cap L^1(\R^d;\R)$ satisfying
\begin{equation}\label{eq:radial-sandwich}
M_\delta^-f\leq f\leq M_\delta^+f,
\qquad
\operatorname{supp}\widehat{M_\delta^\pm f}\subseteq\overline B_2(0,\delta).
\end{equation}
Their one-sided errors satisfy
\begin{equation}\label{eq:width}
\|M_\delta^+f-f\|_1+\|f-M_\delta^-f\|_1
=\omega_{d-1}\int_{]0,+\infty[} u_\nu(\rho,\delta)\,\mathrm d|D\varphi|(\rho)
=W_\delta,
\end{equation}
where $W_\delta$ is the weighted variation cost~\eqref{eq:W}.
\end{proposition}

A proof is given in \ref{app:layers}.

\medskip

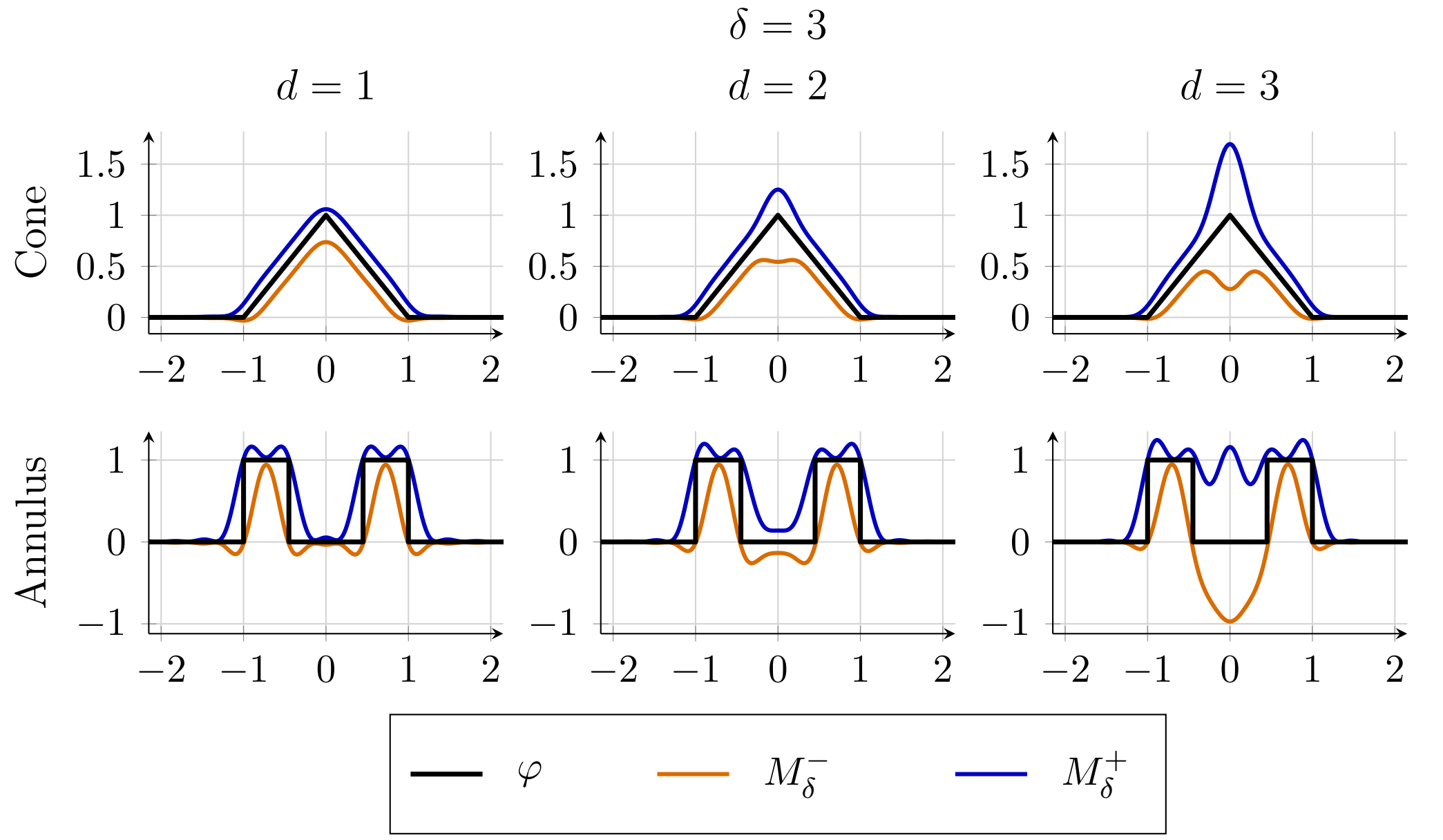
\begin{figure}[htbp]
\centering
\resizebox{\linewidth}{!}{\input{figures/radial_profiles.tex}}
\caption{Radial sections of the Holt--Vaaler-based minorant
and majorant for a common bandwidth $\delta=3$, in dimensions
$d=1,2,3$ (left to right). Top: the conical taper
$\varphi(r)=(1-r)_+$. Bottom: the annular profile
$\varphi(r)=\chi_1(r)-\chi_{0.45}(r)$, with the
normalized half-values at its jumps. Black denotes the target,
orange the constructed minorant $M_\delta^-$, and blue the
constructed majorant $M_\delta^+$. In each dimension the
curves are computed from the radial Holt--Vaaler ball pair through the signed layer formulas~\eqref{eq:envelopes}. The figure-generation code is available in~\cite{GabetRadialVandermondeFigures2026}.}
\label{fig:radial-profiles}
\end{figure}

The signed radial layer construction can be visualized for
both a nonincreasing conical profile and a nonmonotone annulus.
For the conical taper $\varphi(r)=(1-r)_+$, the radial variation
measure is Lebesgue measure on $]0,1[$, so the envelopes in
\eqref{eq:envelopes} are the corresponding integrals of the
Holt--Vaaler ball pair. For the annular profile
$\chi_1-\chi_{0.45}$, the same construction uses the difference
of the ball approximants at the two radii.
Figure~\ref{fig:radial-profiles} displays radial sections of
these envelopes for $d=1,2,3$, using the same Fourier
bandwidth $\delta=3$ in every panel.

\begin{remark}[The layer pair is not generally optimal]
The signed layer construction \eqref{eq:envelopes} is not claimed
to minimize the integrated gap for a general radial profile.
This limitation can already be illustrated by the indicator
$\chi_R$ of a single ball of radius $R>0$. By Holt--Vaaler's
theorem~\cite[Theorem~3]{holtvaaler1996}, the infimum of the
integrated gap over all admissible bandlimited majorant--minorant
pairs equals $\omega_{d-1}u_\nu(R,\delta)$ if and only if
\[
J_\nu(\pi\delta R)J_{\nu+1}(\pi\delta R)=0;
\]
see also \cite[Eq.~(54)]{vaaler2023lattice}.
Consequently, whenever this product is nonzero, there exists an admissible bandlimited majorant--minorant pair for the ball whose integrated gap is strictly smaller than that of our layer pair.
This illustrates why the cost of the particular construction
\eqref{eq:envelopes} is not asserted to be optimal for general
radial profiles.
\end{remark}

For an absolutely continuous profile, \eqref{eq:width} becomes
\[
W_\delta=\omega_{d-1}\int_0^{+\infty}|\varphi'(\rho)|u_\nu(\rho,\delta)\,\mathrm d\rho.
\]
For a decreasing profile, only $\mu^+$ occurs. Both cases use the same construction and the same notation; discontinuous profiles contribute atoms rather than additional boundary terms.

\section{Sampling and proof of the spectral bound}\label{sec:spectral}
\subsection{Exact diagonalization for bandlimited envelopes}
The following identity isolates the role of separation. We state it for general bandlimited functions, since radiality is used only to construct the envelopes.
Related uses of Poisson summation with
ball minorants can be found in
\cite[Eq.~(3.7)]{hockmann2023ShortCommunicationa}.

\begin{lemma}[Bandlimited sampling identity]\label{lem:sampling}
Let $\delta\in \,]0,\Delta]$ be real, and let $H\in L^1(\R^d;\C)$ have Fourier support in $\overline B_2(0,\delta)$, using its continuous Fourier-inversion representative. Then the samples of $H$ on $T^{-1}\mathbb Z^d$ are absolutely summable, and for every $\bm x\in\C^r$,
\begin{equation}\label{eq:diagonal}
\sum_{\bm k\in\mathbb Z^d}H(\bm k/T)
\left|\sum_{j=1}^r x_j\e^{2\pi\ii\langle\bm k/T,\bm\tau_j\rangle}\right|^2
=T^d\left(\int_{\R^d}H\right)\|\bm x\|_2^2.
\end{equation}
\end{lemma}

A proof is given in \ref{app:sampling}.
\medskip

\subsection{Proof of Theorem~\ref{th:main}}
\begin{proof}
\emph{Step 1: Compactly supported profiles.}
Suppose first that $\varphi$ has compact support. Apply
Proposition~\ref{prop:radial} with $\delta=\Delta$ to $f$, and fix
an arbitrary $\bm x\in\C^r$. Multiplying the pointwise sandwich
by the nonnegative squared exponential sum and summing over the
lattice gives
\begin{equation}\label{eq:qf-sandwich}
T^d\left(\int_{\R^d} M_\Delta^-f\right)\|\bm x\|_2^2
\leq\|\bm\Phi_{\bm\tau}\bm x\|_2^2
\leq T^d\left(\int_{\R^d} M_\Delta^+f\right)\|\bm x\|_2^2,
\end{equation}
where Lemma~\ref{lem:sampling} evaluates the envelope sums,
including at $\delta=\Delta$. Their absolute convergence
permits these comparisons even though the lower envelope need
not be nonnegative.

Set
\[
e_\Delta^-:=E-\int_{\R^d} M_\Delta^-f,
\qquad e_\Delta^+:=\int_{\R^d} M_\Delta^+f-E.
\]
Both errors are nonnegative, and \eqref{eq:width} gives
$e_\Delta^-+e_\Delta^+=W_\Delta$. Consequently,
\eqref{eq:qf-sandwich} proves the slightly more informative
enclosure
\begin{equation}\label{eq:asymmetric}
T^d(E-e_\Delta^-)\bm I_r\preceq\bm\Phi_{\bm\tau}^\ast\bm\Phi_{\bm\tau}
\preceq T^d(E+e_\Delta^+)\bm I_r.
\end{equation}
In particular, the asymmetric enclosure \eqref{eq:asymmetric}
implies \eqref{eq:gram} for compactly supported profiles.

\medskip
\noindent\emph{Step 2: Removal of the compact-support assumption.}
Let $\varphi$ now satisfy the general hypotheses of the theorem.
Since $\varphi$ has finite total variation on $]0,+\infty[$,
its limit at infinity exists and is finite. By nonnegativity and
$E<+\infty$, this limit must be zero. Choose a strictly
increasing sequence of continuity points $R_n$
of $\varphi$ with $R_n\to+\infty$, and define
\begin{equation}\label{eq:truncation}
\varphi_n(\rho):=
\begin{cases}
\varphi(\rho),&0\leq\rho<R_n,\\
\frac12\varphi(R_n),&\rho=R_n,\\
0,&\rho>R_n.
\end{cases}
\end{equation}
Each $\varphi_n$ is nonnegative, compactly supported,
and satisfies the representative convention
\eqref{eq:representative}. Furthermore,
$0\leq\varphi_n(\rho)\uparrow\varphi(\rho)$ for every
$\rho\geq0$. Set
\[
E_n:=\omega_{d-1}\int_0^{+\infty}
\varphi_n(\rho)\rho^{d-1}\,\mathrm d\rho,
\qquad
W_{\Delta,n}:=\omega_{d-1}\int_{]0,+\infty[}
u_\nu(\rho,\Delta)\,\mathrm d|D\varphi_n|(\rho).
\]
By monotone convergence, $E_n\uparrow E$ (and $E_n>0$
for all sufficiently large $n$). Because $R_n$ is a
continuity point, the derivative measure satisfies
\[
D\varphi_n=(D\varphi)|_{]0,R_n[}
-\varphi(R_n)\delta_{R_n},
\]
where $(D\varphi)|_{]0,R_n[}$ denotes the restriction of $D\varphi$ to $]0,R_n[$ and $\delta_{R_n}$ denotes the Dirac mass at $R_n$. These two measures are mutually singular, so
\begin{equation}\label{eq:truncation-cost}
W_{\Delta,n}=\omega_{d-1}\left(
\int_{]0,R_n[}u_\nu(\rho,\Delta)
\,\mathrm d|D\varphi|(\rho)
+\varphi(R_n)u_\nu(R_n,\Delta)\right).
\end{equation}
We next verify that the extra boundary term vanishes. Since
$\varphi(+\infty)=0$ and $R_n$ is a continuity point,
\[
0\leq\varphi(R_n)
\leq|D\varphi|(]R_n,+\infty[).
\]
The large-radius asymptotic
$u_\nu(\rho,\Delta)\sim\Delta^{-1}\rho^{d-1}$
from \ref{app:cost}, valid also for $d=1$, implies that
for some constant $C>0$ and all sufficiently large $n$,
$u_\nu(R_n,\Delta)\leq C u_\nu(\rho,\Delta)$
for every $\rho>R_n$. Hence
\[
0\leq\varphi(R_n)u_\nu(R_n,\Delta)
\leq C\int_{]R_n,+\infty[}
 u_\nu(\rho,\Delta)\,\mathrm d|D\varphi|(\rho)
\longrightarrow0,
\]
since the weighted variation in \eqref{eq:W} is finite.
It follows from \eqref{eq:truncation-cost} that
$W_{\Delta,n}\to W_\Delta$.

Let $\bm\Phi_{\bm\tau,n}$ denote the operator given by \eqref{eq:Phi}
with $\varphi_n$ in place of $\varphi$. Step~1 applies to
$\varphi_n$ and gives, for each $\bm x\in\C^r$ and all
sufficiently large $n$,
\[
T^d(E_n-W_{\Delta,n})\|\bm x\|_2^2
\leq\|\bm\Phi_{\bm\tau,n}\bm x\|_2^2
\leq T^d(E_n+W_{\Delta,n})\|\bm x\|_2^2.
\]
The middle term is the nonnegative series
\[
\sum_{\bm k\in\mathbb Z^d}
\varphi_n(\|\bm k/T\|_2)
\left|\sum_{j=1}^r x_j
\e^{2\pi\ii\langle\bm k/T,\bm\tau_j\rangle}
\right|^2.
\]
Its summands increase pointwise to those defined by
$\varphi$. Monotone convergence for series, together
with the boundedness of $E_n+W_{\Delta,n}$, proves
that the limit is finite for every $\bm x\in\C^r$.
Thus \eqref{eq:Phi} defines a linear map
$\bm\Phi_{\bm\tau}:\C^r\to\ell^2(\mathbb Z^d;\C)$,
which is bounded since its domain is finite dimensional.
Passing to the limit in the preceding quadratic-form
inequalities yields \eqref{eq:gram}.

Finally, dividing \eqref{eq:gram} by $T^dE$ proves
\eqref{eq:normalized}. If $W_\Delta<E$, the lower
bound is positive, and taking square roots of the extreme
Gram eigenvalues proves \eqref{eq:kappa}.
\end{proof}

\begin{remark}[The alias restriction]
In dimensions larger than four, the periodic Euclidean distance between two points can exceed $T$, since its maximal value is $T\sqrt d/2$. Thus replacing $\Delta$ in \eqref{eq:sep} by the pairwise separation alone is not justified by this proof. The shortest nonzero diagonal alias always has length $T$, independently of the locations and of the dimension. Even for a single node, the alias-free sampling identity requires the bandwidth restriction $\delta\leq T$.
The same argument also works with the periodic $\ell^\infty$
separation, capped at $T$, because
$\|\bm v\|_2\geq\|\bm v\|_\infty$ for every $\bm v\in\R^d$.
We use the Euclidean separation because it is adapted to the
Euclidean-ball Fourier support of the Holt--Vaaler approximants
and permits a larger admissible bandwidth.
\end{remark}

\begin{remark}[Meaning of the constant]
For compactly supported profiles, $W_\Delta$ is the exact integrated gap of the constructed radial approximation pair in \eqref{eq:envelopes}, and the spectral estimate bounds each one-sided error by this gap. For noncompact profiles, $W_\Delta$ is the limit of the corresponding gaps for the truncations in \eqref{eq:truncation}; no global approximation pair is asserted to exist under the hypotheses of Theorem~\ref{th:main}. Neither argument establishes an optimal spectral constant. The one-dimensional improvement is given in Proposition~\ref{prop:one-dimensional}.
\end{remark}

\begin{remark}[Sampling at discontinuities]\label{rem:jumps}
For compactly supported profiles, the same envelopes bound any representative whose value at a jump located at a radius $\rho>0$ lies between its two traces.
Indeed, approach the jump sphere from both sides in \eqref{eq:radial-sandwich} and use continuity of the envelopes. The integrals and continuous energy are unchanged. Theorem~\ref{th:main} remains valid for these alternative jump values also in the noncompact case: choose the continuity radii in \eqref{eq:truncation}, preserve the prescribed values at all earlier jumps, and repeat the monotone-convergence argument. An arbitrary isolated change at the origin is different: it is not encoded by $D\varphi$ on $]0,+\infty[$. The condition $\varphi(0)=\varphi(0^+)$ cannot simply be omitted.
\end{remark}

\subsection{A sharper one-dimensional specialization}\label{sec:one-dimensional-proof}

\begin{proof}[Proof of Proposition~\ref{prop:one-dimensional}]
Set
\[
f(x):=\varphi(|x|),
\qquad x\in\R.
\]
By the assumptions on $\varphi$ and the representative convention
\eqref{eq:representative}, the function $f$ belongs to
$L^1(\R;\R)$ and is of bounded variation. Moreover,
\[
\int_\R f(x)\,\mathrm dx=E,
\qquad
V_f=2|D\varphi|(]0,+\infty[),
\]
where $V_f$ denotes the total variation
of $f$ on $\R$.

First assume that $\varphi$ has compact support.
Take $\delta=\Delta>0$.
Applying Vaaler's Corollary~12
\cite[Eqs.~(4.8)--(4.9),
p.~201]{vaaler1985}
to $f$ and $-f$, respectively, gives
integrable entire functions
$f_\delta^+$ and $f_\delta^-$ of
exponential type at most $2\pi\delta$
such that
\[
f_\delta^-(x)\leq f(x)\leq f_\delta^+(x),
\qquad x\in\R,
\]
and
\begin{equation}\label{eq:one-dimensional-BV-errors}
\begin{aligned}
\int_\R(f_\delta^+-f)\,\mathrm dx
&=\frac{V_f}{2\delta},\\
\int_\R(f-f_\delta^-)\,\mathrm dx
&=\frac{V_f}{2\delta}.
\end{aligned}
\end{equation}
The minorant here is obtained by taking the negative of Vaaler's majorant of $-f$.

Since these functions are integrable and of exponential type at most $2\pi\delta$, their Fourier transforms are supported in $[-\delta,\delta]$.

By \eqref{eq:one-dimensional-BV-errors},
\[
\int_\R f_\delta^\pm(x)\,\mathrm dx
=
E\pm
\frac{|D\varphi|(]0,+\infty[)}{\delta}.
\]

For an arbitrary $\bm x\in\C^r$, we multiply the pointwise sandwich by
\[
\left|
\sum_{j=1}^r
x_j\e^{2\pi\ii(k/T)\tau_j}
\right|^2
\]
and sum over $k\in\mathbb Z$.
Applying Lemma~\ref{lem:sampling} to
$f_\delta^-$ and $f_\delta^+$ gives
\[
\begin{aligned}
T\left(
E-\frac{|D\varphi|(]0,+\infty[)}{\delta}
\right)\|\bm x\|_2^2
&\leq
\|\bm\Phi_{\bm\tau}\bm x\|_2^2\\
&\leq
T\left(
E+\frac{|D\varphi|(]0,+\infty[)}{\delta}
\right)\|\bm x\|_2^2.
\end{aligned}
\]

Since this holds for every $\bm x\in\C^r$, we obtain
\[
\left\|
\frac{\bm\Phi_{\bm\tau}^\ast\bm\Phi_{\bm\tau}}{TE}
-\bm I_r
\right\|
\leq
\frac{|D\varphi|(]0,+\infty[)}
{E\delta}.
\]
This establishes \eqref{eq:one-dimensional} in the
compactly supported case.

For a general admissible $\varphi$, use the truncations
$\varphi_n$ from \eqref{eq:truncation}. The identity for
$D\varphi_n$ in the proof of Theorem~\ref{th:main}
implies
\[
|D\varphi_n|(]0,+\infty[)
=|D\varphi|(]0,R_n[)+\varphi(R_n)
\longrightarrow |D\varphi|(]0,+\infty[).
\]
Moreover, $E_n\to E>0$ and, for every $\bm x\in\C^r$,
$\|\bm\Phi_{\bm\tau,n}\bm x\|_2^2\to
\|\bm\Phi_{\bm\tau}\bm x\|_2^2$ by monotone convergence,
as shown in the same proof. Passing to the limit in the
compact-support quadratic-form inequalities establishes
\eqref{eq:one-dimensional} without a compact-support
assumption.

\end{proof}

\begin{remark}[On radial lifting]
Although a radial function is determined by its
restriction to a line through the origin, a
one-dimensional bandlimited approximant cannot
in general be extended radially while preserving
integrability. Indeed, $h\in L^1(\R;\R)$ does not
imply $h(\|\cdot\|_2)\in L^1(\R^d;\R)$ for $d\geq2$,
as the latter requires
\[
\int_0^{+\infty} |h(r)|r^{d-1}\,\mathrm dr<+\infty.
\]
For example, $h(t)=(\sin(\pi t)/(\pi t))^2$ is
integrable and bandlimited on $\R$, but its radial
extension is not integrable in any dimension
$d\geq2$. The radial problem therefore involves
a weighted one-dimensional approximation problem,
as reflected in the approach of Holt and
Vaaler~\cite[Lemmas~18--19]{holtvaaler1996}.
\end{remark}

\section{Discussion and perspectives}\label{sec:discussion}

We have established spectral bounds for weighted multivariate Vandermonde operators with nonnegative radial profiles of bounded variation, finite energy, and finite weighted variation cost, without requiring compact support. The proof combines a signed ball-layer decomposition and the bandlimited approximants of Holt and Vaaler for compactly supported profiles, followed by a truncation argument for the general case. The resulting bounds hold in arbitrary dimension, including for nonmonotone and discontinuous profiles, and depend explicitly on the radial variation and the periodic Euclidean separation. The bandlimited sampling identity eliminates both off-diagonal interactions and nonzero diagonal aliases. We have also obtained an elementary upper bound for the Holt--Vaaler cost, exhibiting two regimes, and a sharper one-dimensional specialization.

Several directions remain open. Sharper estimates of the individual upper and lower approximation errors could improve the spectral bounds obtained from their integrated gap. The one-dimensional specialization illustrates the potential benefit of controlling these errors separately, although it relies on a different approximation construction. %
Treating nonradial profiles would require bandlimited approximants adapted to their
geometry.

\clearpage
\appendix
\numberwithin{equation}{section}
\numberwithin{figure}{section}
\numberwithin{table}{section}
\renewcommand{\theequation}{\Alph{section}.\arabic{equation}}
\renewcommand{\thefigure}{\Alph{section}.\arabic{figure}}
\renewcommand{\thetable}{\Alph{section}.\arabic{table}}
\renewcommand{\thesubsection}{\Alph{section}.\arabic{subsection}}

\section{Pointwise radial envelopes and integrability}\label{app:layers}
We give the details behind Proposition~\ref{prop:radial}. The derivative measure is finite, and its support is contained in $]0,L]$ for some $L$. Since the profile vanishes beyond $L$, the Stieltjes identity is
\[
\mu(]s,+\infty[)=\varphi(s^+),\qquad
\mu(\{s\})=\varphi(s^-)-\varphi(s^+)\quad(s>0).
\]
It gives \eqref{eq:layers}, including at zero by continuity from above of the finite signed measure. This also explains why the derivative is taken on the open half-line.

Let $K_\delta=G_\delta-F_\delta\geq0$. From \eqref{eq:ball-pair},
\[
|F_\delta(\bm\xi;\rho)|,\ |G_\delta(\bm\xi;\rho)|
\leq\chi_\rho(\bm\xi)+K_\delta(\bm\xi;\rho).
\]
Consequently,
\begin{equation}\label{eq:L1-control}
\|F_\delta(\cdot;\rho)\|_1,\ \|G_\delta(\cdot;\rho)\|_1
\leq\frac{\omega_{d-1}}d\rho^d+\omega_{d-1}u_\nu(\rho,\delta).
\end{equation}
The right-hand side is bounded on $0<\rho\leq L$. The measurable explicit family in Lemma~\ref{lem:ball-pair}, viewed in the separable space $L^1(\R^d;\R)$, is therefore integrable against each finite Jordan measure. Alternatively, Fubini's theorem applied to \eqref{eq:L1-control} directly constructs these $L^1(\R^d;\R)$ integrals.

There is also a uniform pointwise bound. If $H$ is any continuous integrable function with Fourier support in $\overline B_2(0,\delta)$, Fourier inversion gives
\begin{equation*}%
\|H\|_\infty\leq |B_2(0,\delta)|\,\|\widehat H\|_\infty
\leq |B_2(0,\delta)|\,\|H\|_1.
\end{equation*}
Apply this to the two ball functions and use \eqref{eq:L1-control}. Dominated convergence now shows that the pointwise integrals in \eqref{eq:envelopes} are continuous and agree with their $L^1(\R^d;\R)$ versions. The same bound applied to differences, or Fourier inversion, also identifies their continuous representatives unambiguously.

The inequalities follow by integrating the ball sandwich against $\mu^+$ and $\mu^-$ with the signs in \eqref{eq:envelopes}. Fubini permits Fourier transformation under these integrals, and the Fourier transforms vanish outside the common ball. Finally,
\[
M_\delta^+f-M_\delta^-f
=\int K_\delta(\cdot;\rho)\,\mathrm d|D\varphi|(\rho)\geq0.
\]
The sandwich identifies the sum of the one-sided $L^1(\R^d;\R)$ errors with the integral of this difference. Tonelli and \eqref{eq:ball-gap} prove \eqref{eq:width}.

\section{Absolute sampling and Poisson summation}\label{app:sampling}

We apply the Poisson summation formula to the bandlimited function $H\in L^1(\R^d;\C)$. To justify its pointwise use, we first establish absolute summability of the lattice samples and then identify the Fourier series of the periodized transform. Choose a function in $C_c^\infty(\R^d;\R)$ equal to one on a neighborhood of $\operatorname{supp}\widehat H$, and let $\eta\in\mathcal S(\R^d;\C)$ be its inverse Fourier transform. Then $\eta$ is Schwartz and $H=H*\eta$, pointwise for the continuous representative. Therefore
\begin{equation}\label{eq:sample-sum}
\begin{aligned}
\sum_{\bm k\in\mathbb Z^d}|H(\bm k/T)|
&\leq\int_{\R^d}|H(\bm y)|\sum_{\bm k\in\mathbb Z^d}|\eta(\bm k/T-\bm y)|\,\mathrm d\bm y\\
&\leq C_{\eta,T}\|H\|_1<+\infty.
\end{aligned}
\end{equation}
The inner sum is uniformly bounded: it is periodic in $\bm y$ with periods $T^{-1}\mathbb Z^d$, and Schwartz decay makes the sum uniformly convergent on a fundamental cube.

Now define the continuous $T$-periodic function
\[
P(\bm t):=\sum_{\bm q\in\mathbb Z^d}\widehat H(\bm t+T\bm q),\qquad\bm t\in\R^d.
\]
The sum is locally finite because $\widehat H$ has compact support. For every Fourier index $\bm k\in\mathbb Z^d$, with the Fourier-series basis $\e^{-2\pi\ii\langle\bm k/T,\bm t\rangle}$, its coefficient at $\bm k$ is
\begin{equation*}%
\begin{aligned}
\frac1{T^d}\int_{[0,T]^d}P(\bm t)\e^{2\pi\ii\langle\bm k/T,\bm t\rangle}\,\mathrm d\bm t
&=\frac1{T^d}\int_{\R^d}\widehat H(\bm t)
\e^{2\pi\ii\langle\bm k/T,\bm t\rangle}\,\mathrm d\bm t
=T^{-d}H(\bm k/T).
\end{aligned}
\end{equation*}
The series with these coefficients converges absolutely and uniformly by \eqref{eq:sample-sum}. Uniqueness of Fourier coefficients for continuous periodic functions identifies it with $P$. This proves the required Poisson identity at every $\bm t$.

Because $H\in L^1(\R^d;\C)$, its Fourier transform is continuous. The support assumption implies $\widehat H(\bm t)=0$ for $\|\bm t\|_2>\delta$; continuity then gives the same vanishing on $\|\bm t\|_2=\delta$. Thus $\widehat H(\bm t)=0$ whenever $\|\bm t\|_2\geq\delta$.

Expand the square in \eqref{eq:diagonal}. Each term is evaluated at $\bm t=\bm\tau_\ell-\bm\tau_j$. If $j\ne\ell$, all $\bm t+T\bm q$ have length at least $\Delta\geq\delta$, and so their Fourier terms vanish. If $j=\ell$ and $\bm q\ne0$, their length is at least $T\geq\Delta\geq\delta$, and these terms vanish as well. The only surviving contribution is $\widehat H(\bm 0)=\int_{\R^d}H(\bm\xi)\,\mathrm d\bm\xi$ for each $j=\ell$, proving the lemma, including when $\delta=\Delta$. This argument applies equally to signed envelopes.

\section{Evaluation and limiting regimes of the Holt--Vaaler cost}\label{app:cost}

The fundamental formulas for the Holt--Vaaler cost are recalled in \cite[Section~3]{vaaler2023lattice}. We collect here these known identities and derive their limiting regimes for use in the proof of Theorem~\ref{th:main}.

For fixed $d\in\N^\ast$ and real parameters $\rho,\delta\in \,]0,+\infty[$, set $z:=\pi\delta\rho$ and $\nu:=(d-2)/2$. The equivalent finite expression from \cite[equations~(36)--(40)]{vaaler2023lattice} is
\begin{equation}\label{eq:finite-u}
u_\nu(\rho,\delta)=
\frac{2\delta^{-1}\rho^{d-1}}
 {\pi\{z[J_\nu(z)^2+J_{\nu+1}(z)^2]-(d-1)J_\nu(z)J_{\nu+1}(z)\}}.
\end{equation}
Positivity and continuity follow from the Holt--Vaaler cost theorem. The formula also gives the scaling law
\begin{equation}\label{eq:scaling}
u_\nu(\rho,\delta)=\delta^{-d}u_\nu(\delta\rho,1).
\end{equation}
The origin value is
\begin{equation}\label{eq:origin}
u_\nu(0,\delta)=\frac{d\,2^{d-1}\Gamma(d/2)^2}{\pi^d}\delta^{-d}.
\end{equation}
These facts also ensure that $W_\delta$ is finite for every fixed $\delta>0$ when the radial profile has compact support and finite variation.

For completeness, we detail the small-radius behavior using the differential identity from \cite[equation~(38)]{vaaler2023lattice}. For each $z\in \,]0,+\infty[$ set
\[
A(z)=z[J_\nu(z)^2+J_{\nu+1}(z)^2]-(d-1)J_\nu(z)J_{\nu+1}(z).
\]
The recurrences $J_\nu'(z)=\nu J_\nu(z)/z-J_{\nu+1}(z)$ and $J_{\nu+1}'(z)=J_\nu(z)-(\nu+1)J_{\nu+1}(z)/z$ give
\[
A'(z)=(d-1)J_\nu(z)J_{\nu+1}(z)/z.
\]
If the integer $d>1$, the small-argument expansion gives $\lim_{z\to0^+}A(z)=0$, so \eqref{eq:finite-u} becomes
\begin{equation}\label{eq:finite-integral}
u_\nu(\rho,\delta)=
\frac{\delta^{-1}\rho^{d-1}}
 {\frac{\pi(d-1)}2\displaystyle\int_0^{\pi\delta\rho}
 J_\nu(t)J_{\nu+1}(t)\,\frac{\mathrm dt}{t}}.
\end{equation}
The leading product is
\[
J_\nu(t)J_{\nu+1}(t)
=\frac{t^{d-1}}{2^{d-1}\Gamma(d/2)\Gamma(d/2+1)}\bigl(1+O_d(t^2)\bigr),\qquad t\to0^+.
\]
Substituting it into \eqref{eq:finite-integral} yields \eqref{eq:origin} and
\[
u_\nu(\rho,\delta)=u_\nu(0,\delta)\,\bigl[1+O_d\bigl((\delta\rho)^2\bigr)\bigr],\qquad\delta\rho\to0^+.
\]
In both $O_d$ terms, $d$ is fixed. The finite-integral representation is not used for $d=1$, where its prefactor vanishes.

For $d=1$, the half-integer Bessel formulas give $A(z)=2/\pi$, hence $u_{-1/2}=\delta^{-1}$. For $d=3$, the half-integer Bessel formulas give, as in \cite[Corollary~1.2]{vaaler2023lattice},
\[
u_{1/2}(\rho,\delta)=\frac{\delta^{-1}\rho^2}{1-[\sin(\pi\delta\rho)/(\pi\delta\rho)]^2}.
\]
For fixed $d\in\N^\ast$, the known large-argument Bessel asymptotics, as used in \cite[equation~(37)]{vaaler2023lattice}, give $A(z)\to2/\pi$ as $z\to+\infty$; hence $u_\nu(\rho,\delta)\sim\delta^{-1}\rho^{d-1}$ as $\delta\rho\to+\infty$. The small-radius regime has order $\delta^{-d}$, while the large-radius regime has the order of a boundary layer. Proposition~\ref{prop:elementary} supplies a global elementary majorant, including the transition between these two regimes.

\input{elementary_bound}

\input{comparisons}

\section*{Funding}

This work was supported in part by ANR [grant number ANR-24-CE48-3094],
by AID [grant number AID-2023639], and by the Orange chair on ``Sustainable 6G''
held by CentraleSup\'elec. Kiryung Lee was supported in part by the National
Science Foundation (NSF) CAREER Award [grant number CCF-1943201].
The funders had no role in the study design, analysis, interpretation of results or writing of the manuscript.

\bibliographystyle{elsarticle-num}
\bibliography{references}
\end{document}

%% file: figures/radial_lowpass.tex
\begin{tikzpicture}
\begin{loglogaxis}[width=.96\linewidth,height=7.7cm,
font=\normalsize,tick label style={font=\normalsize},label style={font=\normalsize},
xlabel={$R\Delta$},ylabel={Relative spectral error},xmin=0.25,xmax=8,ymin=0.003,ymax=1500,
legend style={font=\normalsize,at={(0.5,1.04)},anchor=south,legend columns=4},
grid=major,log basis x=2,xtick={0.25,0.5,1,2,4,8},xticklabels={$0.25$,$0.5$,$1$,$2$,$4$,$8$}]
\addplot[thick,blue!75!black,no marks] table[x=RDelta,y=bound] {reproducibility/data/lowpass_bound_d1.dat};
\addlegendentry{$d=1$}
\addplot[only marks,mark=*,mark size=1.2pt,blue!75!black,forget plot] table[x=RDelta,y=deviation] {reproducibility/data/lowpass_empirical_d1.dat};
\addplot[thick,orange!85!black,no marks] table[x=RDelta,y=bound] {reproducibility/data/lowpass_bound_d2.dat};
\addlegendentry{$d=2$}
\addplot[only marks,mark=*,mark size=1.2pt,orange!85!black,forget plot] table[x=RDelta,y=deviation] {reproducibility/data/lowpass_empirical_d2.dat};
\addplot[thick,green!50!black,no marks] table[x=RDelta,y=bound] {reproducibility/data/lowpass_bound_d3.dat};
\addlegendentry{$d=3$}
\addplot[only marks,mark=*,mark size=1.2pt,green!50!black,forget plot] table[x=RDelta,y=deviation] {reproducibility/data/lowpass_empirical_d3.dat};
\addplot[thick,violet,no marks] table[x=RDelta,y=bound] {reproducibility/data/lowpass_bound_d4.dat};
\addlegendentry{$d=4$}
\addplot[only marks,mark=*,mark size=1.2pt,violet,forget plot] table[x=RDelta,y=deviation] {reproducibility/data/lowpass_empirical_d4.dat};
\addplot[black,densely dotted,line width=1.25pt,domain=0.25:8,forget plot] {1};
\addplot[only marks,mark=triangle*,blue!75!black,forget plot] coordinates {(8,0.0035)};
\node[anchor=south east,font=\normalsize] at (axis cs:7.8,0.0035) {$0$};
\end{loglogaxis}
\end{tikzpicture}

%% file: figures/elementary_cost.tex
\begin{tikzpicture}
\begin{loglogaxis}[
  name=cost,
  width=.96\linewidth,
  height=6.1cm,
  tick label style={font=\small},
  label style={font=\normalsize},
  xmin=0.01,xmax=8,
  ymin=1,ymax=1e7,
  xlabel={$\rho$},
  ylabel={Radial cost},
  grid=major,
  major grid style={draw=black!17},
  legend style={font=\small,at={(0.02,0.98)},anchor=north west,
    legend columns=2,draw=black,fill=white,fill opacity=1,text opacity=1,line width=0.3pt},
]
\addplot[blue!75!black,line width=1.15pt,no marks]
  table[x=rho,y=exact] {reproducibility/data/elementary_d8.dat};
\addlegendentry{Exact cost}
\addplot[orange!85!black,line width=1.15pt,
  dash pattern=on 5pt off 2pt on 1.2pt off 2pt,no marks]
  table[x=rho,y=elementary] {reproducibility/data/elementary_d8.dat};
\addlegendentry{Elementary bound}
\addplot[green!45!black,line width=1.05pt,dashed,no marks]
  table[x=rho,y=asymptotic] {reproducibility/data/elementary_d8.dat};
\addlegendentry{$\rho^7$}
\addplot[black,line width=1.05pt,densely dotted,no marks]
  table[x=rho,y=origin] {reproducibility/data/elementary_d8.dat};
\addlegendentry{$u_3(0,1)$}
\end{loglogaxis}

\begin{semilogxaxis}[
  at={(cost.south west)},anchor=north west,yshift=-1.95cm,
  width=.96\linewidth,
  height=4.7cm,
  tick label style={font=\small},
  label style={font=\normalsize},
  xmin=0.01,xmax=8,
  ymin=0.5,ymax=31,
  xlabel={$\rho$},
  ylabel={Bound / exact cost},
  grid=major,
  major grid style={draw=black!17},
  legend style={font=\small,at={(0.02,0.98)},anchor=north west,
    draw=black,fill=white,fill opacity=1,text opacity=1,line width=0.3pt},
]
\addplot[orange!85!black,line width=1.25pt,no marks]
  table[x=rho,y=ratio] {reproducibility/data/elementary_d8.dat};
\addlegendentry{Bound / exact cost}
\addplot[black,densely dotted,line width=1.05pt]
  coordinates {(2.2281692033,0.5) (2.2281692033,31)};
\addlegendentry{$\rho_\delta=7/\pi$}
\addplot[black,dashed,line width=1.05pt]
  coordinates {(0.01,1) (8,1)};
\addlegendentry{$y=1$}
\end{semilogxaxis}
\end{tikzpicture}

%% file: figures/radial_profiles.tex
\begin{tikzpicture}
\begin{groupplot}[
  group style={group size=3 by 2, horizontal sep=0.95cm, vertical sep=0.95cm},
  width=5.03cm, height=3.55cm,
  xmin=-2.15,xmax=2.15,
  xtick={-2,-1,0,1,2},
  axis lines=left,
  tick label style={font=\small},
  label style={font=\normalsize},
  major grid style={draw=black!17},
  grid=major,
  every axis plot/.append style={no marks,line join=round,line cap=round},
  clip=true,
]
\nextgroupplot[ymin=-0.16,ymax=1.82,ytick={0,0.5,1,1.5}]
\addplot[orange!85!black,line width=1.15pt] table[x=x,y=minor] {figures/radial_profiles/cone_d1.dat};
\addplot[blue!75!black,line width=1.15pt] table[x=x,y=major] {figures/radial_profiles/cone_d1.dat};
\addplot[black,line width=1.35pt] table[x=x,y=phi] {figures/radial_profiles/cone_d1.dat};
\nextgroupplot[ymin=-0.16,ymax=1.82,ytick={0,0.5,1,1.5}]
\addplot[orange!85!black,line width=1.15pt] table[x=x,y=minor] {figures/radial_profiles/cone_d2.dat};
\addplot[blue!75!black,line width=1.15pt] table[x=x,y=major] {figures/radial_profiles/cone_d2.dat};
\addplot[black,line width=1.35pt] table[x=x,y=phi] {figures/radial_profiles/cone_d2.dat};
\nextgroupplot[ymin=-0.16,ymax=1.82,ytick={0,0.5,1,1.5}]
\addplot[orange!85!black,line width=1.15pt] table[x=x,y=minor] {figures/radial_profiles/cone_d3.dat};
\addplot[blue!75!black,line width=1.15pt] table[x=x,y=major] {figures/radial_profiles/cone_d3.dat};
\addplot[black,line width=1.35pt] table[x=x,y=phi] {figures/radial_profiles/cone_d3.dat};
\nextgroupplot[ymin=-1.12,ymax=1.35,ytick={-1,0,1}]
\addplot[orange!85!black,line width=1.15pt] table[x=x,y=minor] {figures/radial_profiles/ring_d1.dat};
\addplot[blue!75!black,line width=1.15pt] table[x=x,y=major] {figures/radial_profiles/ring_d1.dat};
\addplot[black,line width=1.35pt] coordinates {(-2.4,0) (-1,0) (-1,1) (-0.45,1) (-0.45,0) (0.45,0) (0.45,1) (1,1) (1,0) (2.4,0)};
\nextgroupplot[ymin=-1.12,ymax=1.35,ytick={-1,0,1}]
\addplot[orange!85!black,line width=1.15pt] table[x=x,y=minor] {figures/radial_profiles/ring_d2.dat};
\addplot[blue!75!black,line width=1.15pt] table[x=x,y=major] {figures/radial_profiles/ring_d2.dat};
\addplot[black,line width=1.35pt] coordinates {(-2.4,0) (-1,0) (-1,1) (-0.45,1) (-0.45,0) (0.45,0) (0.45,1) (1,1) (1,0) (2.4,0)};
\nextgroupplot[ymin=-1.12,ymax=1.35,ytick={-1,0,1}]
\addplot[orange!85!black,line width=1.15pt] table[x=x,y=minor] {figures/radial_profiles/ring_d3.dat};
\addplot[blue!75!black,line width=1.15pt] table[x=x,y=major] {figures/radial_profiles/ring_d3.dat};
\addplot[black,line width=1.35pt] coordinates {(-2.4,0) (-1,0) (-1,1) (-0.45,1) (-0.45,0) (0.45,0) (0.45,1) (1,1) (1,0) (2.4,0)};
\end{groupplot}
\node[font=\normalsize] at ($(group c1r1.north)+(0,0.46)$) {$d=1$};
\node[font=\normalsize] at ($(group c2r1.north)+(0,0.46)$) {$d=2$};
\node[font=\normalsize] at ($(group c3r1.north)+(0,0.46)$) {$d=3$};
\node[font=\normalsize] at ($(group c2r1.north)+(0,1.05)$) {$\delta=3$};
\node[font=\normalsize,rotate=90] at ($(group c1r1.west)+(-1.15cm,0)$) {Cone};
\node[font=\normalsize,rotate=90] at ($(group c1r2.west)+(-1.15cm,0)$) {Annulus};
\node[draw=black,fill=white,inner sep=5pt,anchor=north]
  at ($(group c2r2.south)+(0,-0.78)$) {%
  \begin{tikzpicture}[baseline=(current bounding box.center)]
    \draw[black,line width=1.35pt] (0,0) -- (0.7,0);
    \node[anchor=west,font=\small] at (0.85,0) {$\varphi$};
    \draw[orange!85!black,line width=1.15pt] (2.4,0) -- (3.1,0);
    \node[anchor=west,font=\small] at (3.25,0) {$M^-_{\delta}$};
    \draw[blue!75!black,line width=1.15pt] (5.3,0) -- (6,0);
    \node[anchor=west,font=\small] at (6.15,0) {$M^+_{\delta}$};
  \end{tikzpicture}%
};
\end{tikzpicture}

%% file: elementary_bound.tex
\section{Proof of the elementary two-regime bound}
\label{app:elementary}

Throughout this section, $d\in\N^\ast$ with $d\geq2$ and $\delta>0$ are fixed. We write
\[
m:=d-1\in\N^\ast,
\qquad
\nu:=\frac{m-1}2=\frac{d-2}2,
\qquad
x:=\pi\delta\rho\geq0,
\]
and, for $x>0$,
\[
A(x):=x\bigl[J_\nu(x)^2+J_{\nu+1}(x)^2\bigr]
-mJ_\nu(x)J_{\nu+1}(x),
\]
which is the function $A$ of \ref{app:cost}. By \eqref{eq:finite-u}, $A(x)>0$ for $x>0$, and the cost is $u_\nu(\rho,\delta)=2\delta^{-1}\rho^m/(\pi A(x))$. We work with the dimensionless normalization
\begin{equation}\label{eq:U-normalization}
U(x):=\pi^m\delta^{m+1}
u_\nu\bigl(x/(\pi\delta),\delta\bigr)
=\frac{2x^m}{\pi A(x)}
\qquad (x>0),
\end{equation}
which removes the dependence on $\delta$, in accordance with the scaling law \eqref{eq:scaling}. By \eqref{eq:origin}, $U$ extends continuously to $x=0$ with
\[
U(0)
=\frac{(m+1)2^m\Gamma((m+1)/2)^2}{\pi}.
\]
In these variables, the threshold $\rho_\delta$ of Proposition~\ref{prop:elementary} corresponds to $x=m$, the constant $c_d\delta^{-d}$ to $U(0)$, and the proposition is equivalent to
\begin{equation}\label{eq:U-target}
U(x)\leq
\begin{cases}
U(0)+\bigl(2m^m-U(0)\bigr)\dfrac xm,
&0\leq x\leq m,\\[3mm]
2x^m,
&x\geq m.
\end{cases}
\end{equation}
Both branches take the value $2m^m$ at $x=m$.

The proof is organized as follows. Section~\ref{app:elementary-energy} proves the second branch of \eqref{eq:U-target}, together with a bound valid for all $x>m/2$. Sections~\ref{app:elementary-origin} and~\ref{app:elementary-slope} prove the first branch on $[0,x_0]$, where
\[
x_0:=\frac{m+1}2\in\,\Bigl]\frac m2,m\Bigr].
\]
Section~\ref{app:elementary-closing} covers the remaining interval $[x_0,m]$.

\subsection{An energy inequality beyond the turning scale}
\label{app:elementary-energy}

For $x>0$, put
\[
p(x):=\sqrt{x}\,J_\nu(x),
\qquad
q(x):=\sqrt{x}\,J_{\nu+1}(x).
\]
Since $2\nu+1=m$, the Bessel recurrences
\cite[Eqs.~(10.6.1)--(10.6.2)]{dlmf} give
\[
p'=\frac{m}{2x}p-q,
\qquad
q'=p-\frac{m}{2x}q,
\qquad
A=p^2+q^2-\frac{m}{x}pq,
\]
and a direct computation yields $A'=mpq/x^2$, in agreement with \ref{app:cost}. For $x>m/2$, let
\[
h(x):=1-\frac{m}{2x}>0,
\qquad\text{so that}\qquad
h'(x)=\frac{m}{2x^2}.
\]
Then
\begin{equation}\label{eq:energy-monotonicity}
\left(\frac{A}{h}\right)'
=\frac{A'h-Ah'}{h^2}
=\frac{m}{x^2h^2}
\left(pq\,h-\frac A2\right)
=-\frac{m\,(p-q)^2}{2x^2h^2}
\leq0,
\end{equation}
where the last equality uses $pqh-A/2=-(p-q)^2/2$, which follows from the expression of $A$ in terms of $p$ and $q$.
Hence $A/h$ is nonincreasing on $]m/2,+\infty[$. The large-argument expansion
\cite[Eq.~(10.17.3)]{dlmf}, at fixed order, gives $p^2+q^2\to2/\pi$ and $pq/x\to0$ as $x\to+\infty$, so that
\[
\lim_{x\to+\infty}\frac{A(x)}{h(x)}
=\frac2\pi.
\]
A nonincreasing function is at least its limit, and therefore
\begin{equation}\label{eq:energy-upper-U}
A(x)\geq
\frac2\pi\left(1-\frac{m}{2x}\right),
\qquad
U(x)\leq\frac{2x^{m+1}}{2x-m}
\qquad(x>m/2).
\end{equation}
If $x\geq m$, then $2x-m\geq x$, and \eqref{eq:energy-upper-U} yields $U(x)\leq2x^m$. This proves the second branch of \eqref{eq:U-target}.

\subsection{A near-origin exponential bound}
\label{app:elementary-origin}

We first prove a Bessel lower bound valid for every real order $\alpha\geq0$:
\begin{equation}\label{eq:J-lower}
J_\alpha(x)\geq
\frac{(x/2)^\alpha}{\Gamma(\alpha+1)}
\exp\left(-\frac{x^2}{2(\alpha+1)}\right),
\qquad 0<x\leq\alpha+1.
\end{equation}

Fix $\alpha\geq0$. By the power series
\cite[Eq.~(10.2.2)]{dlmf}, $J_\alpha(x)>0$ for all sufficiently small $x>0$. On any interval $]0,x_1[$ on which $J_\alpha>0$, the quotient
\[
R(x):=\frac{J_{\alpha+1}(x)}{J_\alpha(x)}
\]
is well defined, and the recurrences
\cite[Eq.~(10.6.2)]{dlmf} give the Riccati equation
\[
R'(x)=1-\frac{2\alpha+1}{x}R(x)+R(x)^2,
\]
while the power series gives
\[
R(x)=\frac{x}{2(\alpha+1)}
+O_\alpha(x^3),
\qquad x\to0^+.
\]
We compare $R$ with
\[
\beta(x):=\frac{x}{\alpha+1},
\qquad 0<x\leq\alpha+1.
\]
Near the origin, $R<\beta$. If $R(x)=\beta(x)$ at some point $x\in\,]0,\min(x_1,\alpha+1)[$, then inserting $R=\beta$ into the Riccati equation gives
\[
R'(x)-\beta'(x)
=-1+\frac{x^2}{(\alpha+1)^2}<0.
\]
At the first such contact point, however, $R-\beta$ would be nonpositive on the left and vanish at the point, which forces $R'-\beta'\geq0$ there. This contradiction shows that
\[
R(x)\leq\frac{x}{\alpha+1}
\qquad\text{for }0<x<\min(x_1,\alpha+1).
\]
On this interval, the same recurrences give
\[
\frac{\mathrm d}{\mathrm dx}
\log\left(\frac{J_\alpha(x)}{x^\alpha}\right)
=-R(x)
\geq-\frac{x}{\alpha+1}.
\]
Integrating from $0$, where $J_\alpha(x)/x^\alpha\to2^{-\alpha}/\Gamma(\alpha+1)$ by
\cite[Eq.~(10.2.2)]{dlmf}, proves \eqref{eq:J-lower} for $0<x<\min(x_1,\alpha+1)$.
Finally, $J_\alpha$ has no zero in $]0,\alpha+1]$: if $x_1\leq\alpha+1$ were the first positive zero, letting $x\to x_1^-$ in the bound just proved would give $J_\alpha(x_1)>0$. Thus \eqref{eq:J-lower} holds on all of $]0,\alpha+1[$, and at $x=\alpha+1$ by continuity.

We now apply \eqref{eq:J-lower} with $\alpha=\nu$ and $\alpha=\nu+1$. Both are admissible for $0<x\leq\nu+1=x_0$, and since $2\nu+1=m$, $\Gamma(\nu+1)=\Gamma((m+1)/2)$, and $\Gamma(\nu+2)=\Gamma((m+3)/2)$, the product of the two bounds reads, for $0<x\leq x_0$,
\[
J_\nu(x)J_{\nu+1}(x)
\geq
\frac{x^m}
{2^m\Gamma((m+1)/2)\Gamma((m+3)/2)}
\exp\left(
-\frac{2(m+2)x^2}{(m+1)(m+3)}
\right).
\]
Recall from \ref{app:cost} that $A(x)\to0$ as $x\to0^+$ and that
\[
A'(x)=\frac{m}{x}J_\nu(x)J_{\nu+1}(x).
\]
Integrating the product bound from $0$ to $x$ and substituting $t=xv$ gives
\[
A(x)\geq\frac{m\,x^m}
{2^m\Gamma((m+1)/2)\Gamma((m+3)/2)}
\int_0^1 v^{m-1}\exp\left(
-\frac{2(m+2)x^2v^2}{(m+1)(m+3)}
\right)\mathrm dv.
\]
The identity $\Gamma((m+3)/2)=\frac{m+1}2\Gamma((m+1)/2)$ shows that the prefactor equals $2x^m/(\pi U(0))$. Moreover, $mv^{m-1}\,\mathrm dv$ is a probability measure on $[0,1]$ with second moment $m/(m+2)$, so Jensen's inequality for the convex function $y\mapsto\exp(-cy)$, with $c=2(m+2)x^2/((m+1)(m+3))$, gives
\[
A(x)\geq\frac{2x^m}{\pi U(0)}
\exp\left(-\frac{2mx^2}{(m+1)(m+3)}\right).
\]
By \eqref{eq:U-normalization}, this is equivalent to
\begin{equation}\label{eq:U-exponential}
U(x)\leq
U(0)\exp\left(b\,x^2\right),
\qquad
b:=\frac{2m}{(m+1)(m+3)},
\qquad 0\leq x\leq x_0,
\end{equation}
where the value at $x=0$ follows by continuity.

\subsection{A scalar bound for the interpolation slope}
\label{app:elementary-slope}

Set
\[
\theta:=bx_0^2
=\frac{m(m+1)}{2(m+3)}.
\]
We claim that
\begin{equation}\label{eq:gamma-slope}
\frac{U(0)\e^\theta}{x_0}
\leq 2m^{m-1}.
\end{equation}

The claim rests on a midpoint Stirling inequality, which we prove for completeness. For $y>0$, the integral representation of the digamma function
\cite[Eq.~(5.9.12)]{dlmf} and the Frullani integral
$\log y=\int_0^{+\infty}(\e^{-t}-\e^{-yt})\,\mathrm dt/t$ combine into
\[
\psi(y+\tfrac12)-\log y
=
\int_0^{+\infty}\e^{-yt}
\left(
\frac1t-\frac1{2\sinh(t/2)}
\right)\mathrm dt>0,
\]
the integrand being positive because $2\sinh(t/2)>t$ for $t>0$. Hence the function
\[
y\longmapsto
\log\Gamma(y+\tfrac12)
-y\log y+y-\tfrac12\log(2\pi)
\]
is increasing on $]0,+\infty[$. By Stirling's expansion
\cite[Eq.~(5.11.1)]{dlmf}, it tends to zero as $y\to+\infty$, so it is negative. This is the midpoint Stirling inequality
\begin{equation}\label{eq:midpoint-stirling}
\Gamma(y+\tfrac12)
<\sqrt{2\pi}\,(y/\e)^y,
\qquad y>0.
\end{equation}
With $y=m/2$, \eqref{eq:midpoint-stirling} gives $\Gamma((m+1)/2)^2<2\pi(m/(2\e))^m$, that is,
\[
\frac{2^m\Gamma((m+1)/2)^2}{\pi}
<2m^m\e^{-m}.
\]
Since $(m+1)/x_0=2$, the left-hand side of \eqref{eq:gamma-slope} equals $2\cdot\frac{2^m\Gamma((m+1)/2)^2}{\pi}\e^\theta<4m^m\e^{-m+\theta}$. The claim \eqref{eq:gamma-slope} therefore follows once we show
\[
2m\,\e^{-(m-\theta)}\leq1,
\qquad\text{where}\qquad
m-\theta=\frac{m(m+5)}{2(m+3)}.
\]
Equivalently, we show that
\[
g(m):=
\frac{m(m+5)}{2(m+3)}
-\log(2m)
\]
is positive for every integer $m\geq1$. Viewing $m$ as a real variable, and using $m(m+5)/(m+3)=m+2-6/(m+3)$,
\[
g'(m)
=\frac12+\frac3{(m+3)^2}-\frac1m>0
\qquad\text{for } m\geq2,
\]
so it suffices to check $g(1)>0$ and $g(2)>0$. We have $g(1)=\frac34-\log2$ and $g(2)=\frac75-\log4$. Both are positive, since the partial sums of the exponential series give
\[
\e^{3/4}\geq1+\frac34+\frac{(3/4)^2}{2}>2,
\qquad
\e^{7/5}\geq\sum_{k=0}^5\frac{(7/5)^k}{k!}>4.
\]
This proves \eqref{eq:gamma-slope}.

We can now prove the first branch of \eqref{eq:U-target} on $[0,x_0]$. The function $x\mapsto\e^{bx^2}$ is convex, hence lies below its secant on $[0,x_0]$:
\[
\e^{bx^2}\leq1+(\e^\theta-1)\frac x{x_0},
\qquad 0\leq x\leq x_0.
\]
Combined with \eqref{eq:U-exponential}, this gives
\[
U(x)\leq
U(0)+\frac{U(0)(\e^\theta-1)}{x_0}\,x,
\qquad 0\leq x\leq x_0.
\]
By \eqref{eq:gamma-slope} and $x_0\leq m$, the slope satisfies
\begin{multline*}
\frac{U(0)(\e^\theta-1)}{x_0}
=\frac{U(0)\e^\theta}{x_0}-\frac{U(0)}{x_0}
\leq2m^{m-1}-\frac{U(0)}{x_0}\\
\leq2m^{m-1}-\frac{U(0)}m
=\frac{2m^m-U(0)}m,
\end{multline*}
which is the slope in \eqref{eq:U-target}. The first branch of \eqref{eq:U-target} is thus proved on $[0,x_0]$.

\subsection{Closing the intermediate interval}
\label{app:elementary-closing}

It remains to prove the first branch of \eqref{eq:U-target} on $[x_0,m]$. If $m=1$, then $x_0=m$ and there is nothing to prove. We treat $m\geq3$ first and $m=2$ separately.

Let $m\geq3$ and $x\in[x_0,m]$, and put $t:=x/m\in[\frac{m+1}{2m},1]$. We claim that
\begin{equation}\label{eq:t-polynomial}
t^m\leq2t-1,
\qquad
\frac{m+1}{2m}\leq t\leq1.
\end{equation}
The function $t\mapsto t^m-2t+1$ is convex, so it suffices to check \eqref{eq:t-polynomial} at the two endpoints. At $t=1$ both sides vanish. At $t=(m+1)/(2m)$, the inequality reads $m\bigl(\frac{m+1}{2m}\bigr)^m\leq1$. Since $(m+1)/(2m)\leq2/3$ for $m\geq3$, and since the sequence $m(2/3)^m$ is nonincreasing for $m\geq3$ (the ratio of consecutive terms is $2(m+1)/(3m)\leq8/9$),
\[
m\left(\frac{m+1}{2m}\right)^m
\leq m\left(\frac23\right)^m
\leq3\left(\frac23\right)^3
=\frac89<1.
\]
This proves \eqref{eq:t-polynomial}. Now $x>m/2$, so \eqref{eq:energy-upper-U} applies, and in terms of $t$ it reads $U(x)\leq2m^mt^{m+1}/(2t-1)$. Using \eqref{eq:t-polynomial},
\[
U(x)
\leq2m^m\,\frac{t^{m+1}}{2t-1}
\leq2m^m\,t
=2m^{m-1}x
\leq U(0)+\bigl(2m^m-U(0)\bigr)\frac xm,
\]
where the last inequality is equivalent to $U(0)(1-x/m)\geq0$. This is the first branch of \eqref{eq:U-target} on $[x_0,m]$ when $m\geq3$.

Finally, let $m=2$, so that $x_0=3/2$, $U(0)=3$, and the first branch of \eqref{eq:U-target} reads $U(x)\leq3+\frac52x$. By the $d=3$ formula of \ref{app:cost},
\[
U(x)=
\frac{x^2}{1-(\sin x/x)^2}.
\]
For $3/2\leq x\leq2$,
\[
1-\left(\frac{\sin x}{x}\right)^2
\geq1-\frac1{x^2}
\geq\frac59>\frac12,
\]
hence
\[
U(x)\leq2x^2
\leq4x
\leq3+\frac52x,
\]
where the last two inequalities use $x\leq2$.

The inequality \eqref{eq:U-target} is now proved for every $x\geq0$ and every integer $m\geq1$. Undoing the normalization \eqref{eq:U-normalization} gives Proposition~\ref{prop:elementary}.

%% file: comparisons.tex
\section{Comparison with related spectral bounds}\label{app:comparisons}

We compare the explicit bound of Theorem~\ref{th:main} with two results designed for particular frequency profiles. Each comparison uses the same nodes, sampled Gram matrix, and energy normalization for both estimates.

\subsection{Ball sampling: comparison with Li}\label{app:comparison-li}

\begin{figure}[htbp]
\centering
\resizebox{\linewidth}{!}{\input{figures/comparison_li}}
\caption{Lower bounds and actual smallest eigenvalue of the same
normalized closed-ball Gram matrix, with $d=2$, $R=1$, $T=2\Delta$,
and nodes $\{0,\Delta\}^2$. The figure-generation code is available
in~\cite{GabetRadialVandermondeFigures2026}.}
\label{fig:comparison-li}
\end{figure}
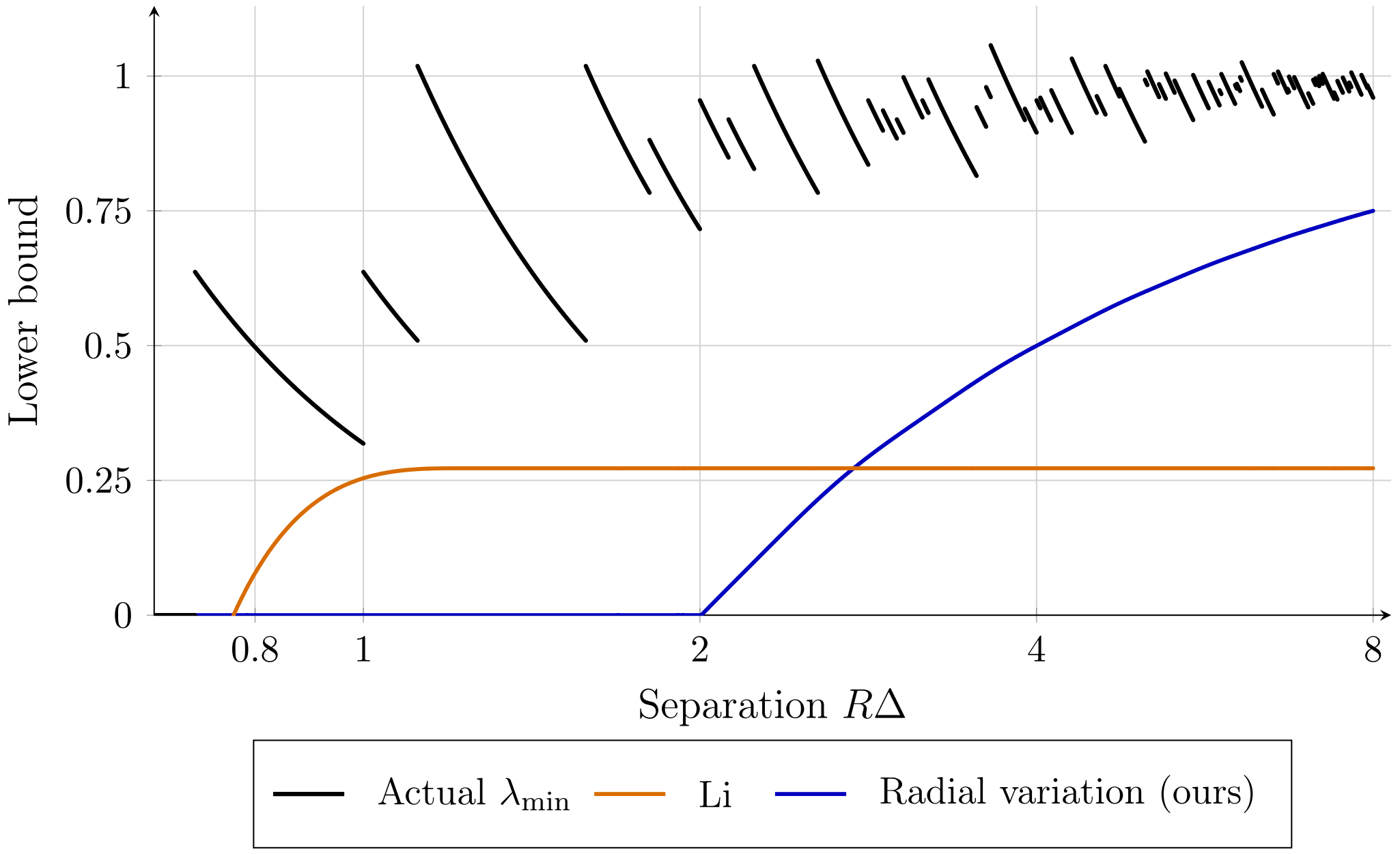

Consider the ideal low-pass profile of radius $R$, with the value $1$
on its cutoff sphere. This closed-ball convention is allowed by
Remark~\ref{rem:jumps}; it is used here to match exactly the discrete
ball-sampling matrix of Li~\cite{li_multivariate2025}. In Li's notation, the frequency
radius equals $RT$ and the nodes are
$\{\bm\tau_j/T\}$. For $d\geq2$,
Li's Theorem~2.5 yields the normalized lower bound
\[
\frac{d\,c_d^{\rm Li}(\alpha)}{\omega_{d-1}},
\qquad
\frac{j_{d/2-1,1}}{\pi}<\alpha\leq R\Delta,
\]
where $j_{\eta,1}$ is the first positive zero of $J_\eta$ and
$c_d^{\rm Li}$ is the constant denoted $c_d$ in~\cite[Eqs.~(2.3)--(2.4)]{li_multivariate2025}.
For this same matrix, Corollary~\ref{cor:ball} gives the
bound compared with Li's in Figure~\ref{fig:comparison-li}. Unlike Li's result, which provides a lower bound on the smallest eigenvalue, our corollary controls both ends of the Gram spectrum.

The estimates are complementary. Li's bound is
positive at smaller separations, whereas the radial
variation bound is larger at greater separations.
The sawtooth behavior of the actual smallest eigenvalue
reflects discrete changes in the lattice frequencies
contained in the closed ball. Li's bound is applicable
only for $R\Delta>j_{0,1}/\pi\approx0.76548$ in this setting.
Moreover, our radial variation bound is asymptotically sharp as
$R\Delta\to+\infty$, since it converges to $1$,
the limiting smallest eigenvalue of the normalized Gram matrix.
By contrast, Li's lower bound remains bounded away from $1$,
even for arbitrarily large separations.

Li also discusses an alternative bound based on the Holt--Vaaler minorant, which yields better constants in high dimensions at the expense of a stronger separation requirement~\cite[discussion following Theorem~2.5]{li_multivariate2025}.
In dimension two, this alternative yields the constant lower bound $1/4$ only for $R\Delta\gtrsim 2.67$, and is therefore weaker than Li's bound displayed in Figure~\ref{fig:comparison-li}.

Related lower bounds based on bandlimited
minorants of Euclidean balls were studied
by Gon\c{c}alves~\cite[Theorem~1]{G2018},
who identified the critical radius for
positive-integral minorants and determined
their optimal integral in a specified range.

\subsection{Gaussian weights: comparison with Gon\c{c}alves--Kelly--Madrid}
\label{app:comparison-gkm}

Theorem~\ref{th:main} also applies directly to the Gaussian
$\varphi(\rho)=\e^{-\pi\rho^2}$, without truncation. In dimension
$d=2$, its continuous energy is $E=1$. The Gaussian approximants
of Gon\c{c}alves, Kelly, and Madrid~\cite{goncalves_kelly_madrid2015}
have Fourier transforms supported in a box. Choosing the half-width
$a=\Delta/\sqrt2$ ensures that every nonzero periodic node difference
or diagonal alias lies outside the interior of this box. The sampling
argument of \ref{app:sampling} therefore applies, with
boundary terms vanishing by continuity of the Fourier transforms.

For the unit-scale Gaussian, set
\[
P(a)=\sum_{n\in\mathbb Z}\e^{-\pi a^2n^2},\qquad
Q(a)=\sum_{n\in\mathbb Z}(-1)^n\e^{-\pi a^2n^2}.
\]
Since $E=1$, Theorem~\ref{th:main} gives the enclosure
$[\max\{0,1-W_\Delta\},1+W_\Delta]$ for the spectrum of the normalized Gram matrix $\bm\Phi_{\bm\tau}^\ast\bm\Phi_{\bm\tau}/T^2$. The integrals
of the box-supported Gaussian majorant and minorant
in~\cite[proofs of Theorems~1--2, Eqs.~(3.6)--(3.11)]{goncalves_kelly_madrid2015}
give the following enclosure for the same normalized sampled Gram matrix:
\[
\operatorname{spec}\!\left(\frac{\bm\Phi_{\bm\tau}^\ast\bm\Phi_{\bm\tau}}{T^2}\right)
\subseteq
\bigl[\max\{0,\,2P(a)Q(a)-P(a)^2\},\;P(a)^2\bigr].
\]

\begin{figure}[htbp]
\centering
\input{figures/comparison_gaussian}
\caption{Eigenvalue enclosures for the same normalized, untruncated
Gaussian Gram matrix with $d=2$, $T=2\Delta$, and nodes
$\{0,\Delta\}^2$. The comparison uses the box half-width
$a=\Delta/\sqrt2$. The figure-generation code is available in~\cite{GabetRadialVandermondeFigures2026}.}
\label{fig:comparison-gaussian}
\end{figure}
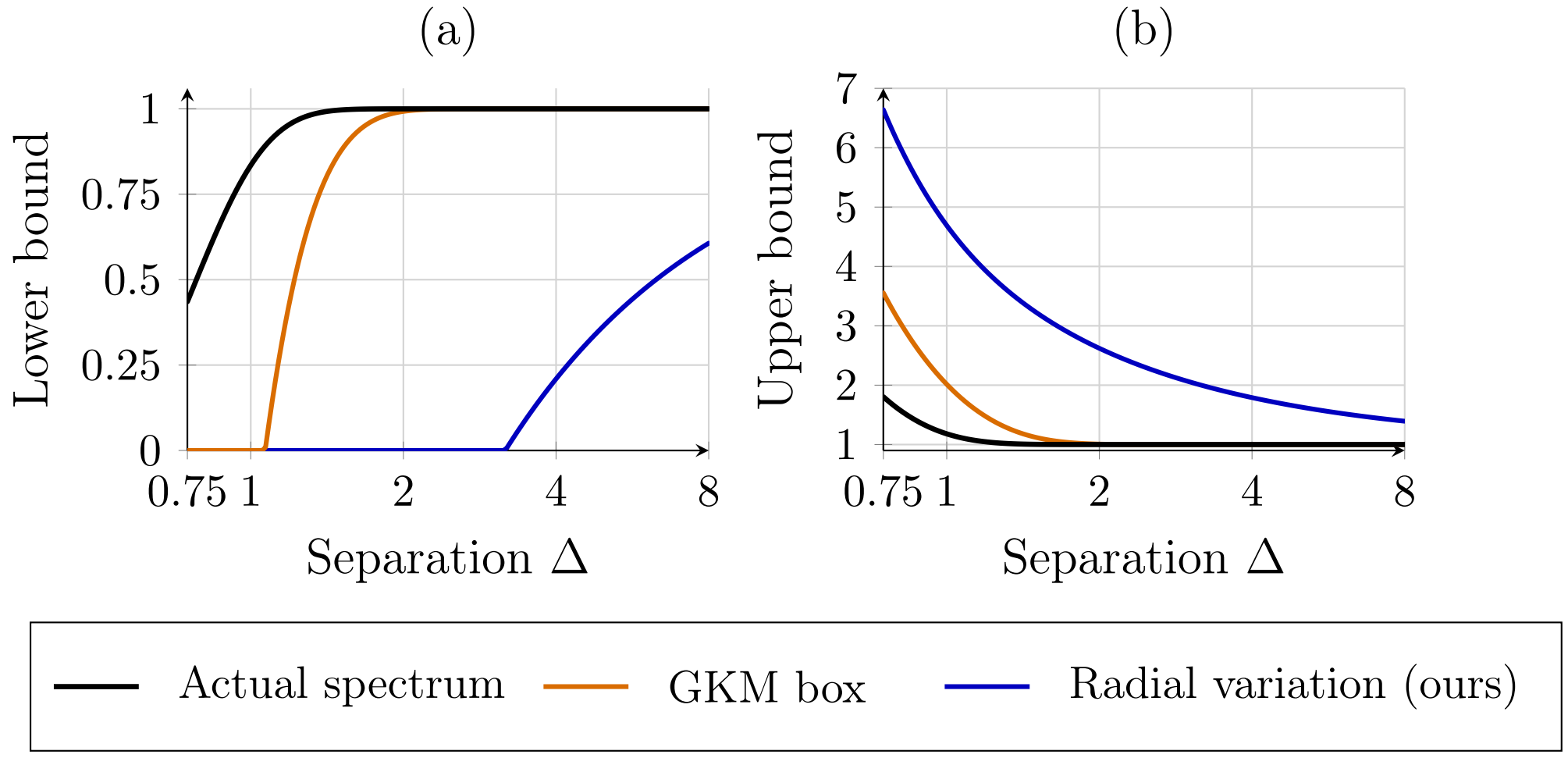

The Gaussian-specific envelopes yield tighter bounds throughout the displayed range. In particular, at $\Delta=2$, our lower bound remains zero even though the sampled Gram matrix is nearly the identity. This illustrates the trade-off between generality and sharpness. The radial variation method applies without a profile-specific extremal construction, but may yield conservative bounds for particular profiles.
Nevertheless, for this two-dimensional Gaussian,
the large-bandwidth asymptotic
$u_0(\rho,\Delta)\sim\rho/\Delta$ and dominated
convergence in \eqref{eq:W} yield
\[
\lim_{\Delta\to+\infty}\Delta W_\Delta
=2\pi\int_0^{+\infty}\rho|\varphi'(\rho)|\,\mathrm d\rho
=\pi.
\]
Consequently, $W_\Delta\sim\pi/\Delta$, so both endpoints
of our eigenvalue enclosure converge to $1$ as
$\Delta\to+\infty$. The error bound therefore decreases
at the rate $1/\Delta$.

%% file: figures/comparison_li.tex
\begin{tikzpicture}
\pgfplotsset{
  comparison axis/.style={
    axis lines=left,
    tick label style={font=\small},
    label style={font=\normalsize},
    major grid style={draw=black!17},
    grid=major,
    every axis plot/.append style={no marks,line join=round,line cap=round},
    unbounded coords=jump,
    clip=true,
  },
  comparison radial/.style={blue!75!black,line width=1.15pt},
  comparison li/.style={orange!85!black,line width=1.15pt},
  comparison exact/.style={black,line width=1.25pt},
}
\begin{axis}[
  name=li-axis,
  comparison axis,
  width=14.1cm,height=7.75cm,
  xmode=log,log basis x=2,
  xmin=0.65,xmax=8.30,
  xtick={0.8,1,2,4,8},
  xticklabels={$0.8$,$1$,$2$,$4$,$8$},
  xlabel={Separation $R\Delta$},
  ymin=0,ymax=1.13,
  ytick={0,0.25,0.5,0.75,1},
  ylabel={Lower bound},
]
\addplot[comparison radial] table[x=s,y=radial_lower] {reproducibility/data/comparison_li.dat};
\addplot[comparison li] table[x=s,y=Li_lower] {reproducibility/data/comparison_li.dat};
\addplot[comparison exact] table[x=s,y=lambda_min] {reproducibility/data/comparison_li.dat};
\end{axis}
\node[draw=black,fill=white,inner sep=5pt,anchor=north]
 at ($(li-axis.south)+(0,-1.26cm)$) {%
 \begin{tikzpicture}[baseline=(current bounding box.center)]
  \draw[black,line width=1.25pt] (0,0) -- (0.72,0);
  \node[anchor=west,font=\small] at (0.87,0) {Actual $\lambda_{\min}$};
  \draw[orange!85!black,line width=1.15pt] (3.25,0) -- (3.97,0);
  \node[anchor=west,font=\small] at (4.12,0) {Li};
  \draw[blue!75!black,line width=1.15pt] (5.08,0) -- (5.80,0);
  \node[anchor=west,font=\small] at (5.95,0) {Radial variation (ours)};
 \end{tikzpicture}%
};
\end{tikzpicture}

%% file: figures/comparison_gaussian.tex
\begin{tikzpicture}
\pgfplotsset{
  gaussian axis/.style={
    axis lines=left,
    tick label style={font=\small},
    label style={font=\normalsize},
    major grid style={draw=black!17},
    grid=major,
    every axis plot/.append style={no marks,line join=round,line cap=round},
    clip=true,
  },
  gaussian radial/.style={blue!75!black,line width=1.15pt},
  gaussian box/.style={orange!85!black,line width=1.15pt},
  gaussian exact/.style={black,line width=1.25pt},
}
\begin{groupplot}[
  group style={group size=2 by 1, horizontal sep=1.48cm},
  width=6.0cm,height=4.65cm,
  gaussian axis,
  xmode=log,log basis x=2,
  xmin=0.75,xmax=8,
  xtick={0.75,1,2,4,8},
  xticklabels={$0.75$,$1$,$2$,$4$,$8$},
  xlabel={Separation $\Delta$},
]
\nextgroupplot[
  ymin=0,ymax=1.06,
  ytick={0,0.25,0.5,0.75,1},
  ylabel={Lower bound},
]
\addplot[gaussian radial] table[x=Delta,y=radial_lower] {reproducibility/data/comparison_gaussian_dense.dat};
\addplot[gaussian box] table[x=Delta,y=gkm_lower] {reproducibility/data/comparison_gaussian_dense.dat};
\addplot[gaussian exact] table[x=Delta,y=actual_lower] {reproducibility/data/comparison_gaussian_dense.dat};

\nextgroupplot[
  ymin=.9,ymax=7,
  ytick={1,2,3,4,5,6,7},
  ylabel={Upper bound},
]
\addplot[gaussian radial] table[x=Delta,y=radial_upper] {reproducibility/data/comparison_gaussian_dense.dat};
\addplot[gaussian box] table[x=Delta,y=gkm_upper] {reproducibility/data/comparison_gaussian_dense.dat};
\addplot[gaussian exact] table[x=Delta,y=actual_upper] {reproducibility/data/comparison_gaussian_dense.dat};
\end{groupplot}

\node[font=\normalsize] at ($(group c1r1.north)+(0,0.48)$) {(a)};
\node[font=\normalsize] at ($(group c2r1.north)+(0,0.48)$) {(b)};

\node[draw=black,fill=white,inner sep=5pt,anchor=north]
  at ($(group c1r1.south)!0.5!(group c2r1.south)+(0,-1.45cm)$) {%
  \begin{tikzpicture}[baseline=(current bounding box.center)]
    \draw[black,line width=1.25pt] (0,0) -- (0.72,0);
    \node[anchor=west,font=\small] at (0.87,0) {Actual spectrum};
    \draw[orange!85!black,line width=1.15pt] (4.15,0) -- (4.87,0);
    \node[anchor=west,font=\small] at (5.02,0) {GKM box};
    \draw[blue!75!black,line width=1.15pt] (7.55,0) -- (8.27,0);
    \node[anchor=west,font=\small] at (8.42,0) {Radial variation (ours)};
  \end{tikzpicture}%
};
\end{tikzpicture}